\input amstex
\documentstyle{amsppt}
\NoBlackBoxes

\def\l{\lambda}
\def\la{\lambda}
\def\ga{\gamma}
\def\U{{\Cal U}}

\def\A{{\Bbb A}}

\def\Fr{{\Bbb F}}
\def\Na{{\Bbb N}}

\def\C{{\Bbb C}}
\def\Al{{\Cal A}}
\def\E{{\Cal E}}
\def\H{{\Cal H}}
\def\L{{\Cal L}}
\def\N{{\Cal N}}
\def\M{{\Cal M}}
\def\F{{\Cal F}}
\def\P{{\Cal P}}
\def\S{{\Cal S}}
\def\W{{\Cal W}}

\def\K{{\Cal K}}
\define\({\left(}
\define\){\right)}
\define\[{\left[}
\define\]{\right]}

\def\longlongrightarrow{\relbar\joinrel\relbar\joinrel
\relbar\joinrel\rightarrow}

\hoffset= .72truein

\topmatter

\title Noncommutative  Interpolation and Poisson Transforms 
\endtitle

\rightheadtext{ Interpolation and  Poisson Transforms }
\leftheadtext{A. Arias and G. Popescu}

\author Alvaro Arias and Gelu Popescu
\endauthor
\thanks  Draft 08/26/97\endthanks

\thanks The second author was partially supported by  NSF
 DMS-9531954.
\endthanks

\address Alvaro Arias
\newline Division of Mathematics and Statistics, 
The University of Texas at San Antonio, San
Antonio, TX 78249, U.S.A.\endaddress
\email arias\@math.utsa.edu\endemail

\address Gelu Popescu
\newline Division of Mathematics and Statistics, 
The University of Texas at San Antonio, San
Antonio, TX 78249, U.S.A.\endaddress
\email gpopescu\@math.utsa.edu\endemail

\abstract
General results of interpolation (eg. Nevanlinna-Pick) by elements in the noncommutative analytic Toeplitz algebra $F^\infty$ 
(resp. noncommutative disc algebra $\Al_n$) with consequences to the interpolation by bounded operator-valued analytic 
functions in the unit ball of $\Bbb C^n$
are obtained.

Non-commutative Poisson transforms are used to provide new von Neumann type inequalities.
Completely isometric representations of the quotient algebra $F^\infty/J$ on Hilbert spaces, where $J$ is any $w^*$-closed, 2-sided ideal of $F^\infty$, are obtained and used to construct a $w^*$-continuous, $F^\infty/J$--functional calculus associated to row contractions $T=[T_1,\dots, T_n]$ when $f(T_1,\dots,T_n)=0$ for any $f\in J$.
Other properties of the dual algebra $F^\infty/J$ are considered.
\endabstract


\endtopmatter

\document

\baselineskip16pt


In [Po5],
the second author proved the following version of von Neumann's inequality 
for row contractions:  if
$T_1,\dots, T_n\in B(\H)$ (the algebra of all bounded  linear operators on the Hilbert space $\H$) and $T=[T_1,\dots,T_n]$ is a contraction, i.e., $\sum_{i=1}^nT_iT_i^*\leq I_\H$,
then for every polynomial $p(X_1,\dots,X_n)$
on $n$ noncommuting indeterminates,
$$\bigl\|p(T_1,\dots,T_n)\bigr\|_{B(\H)}
\leq \bigl\|p(S_1,\dots,S_n)\bigr\|_{B(\F^2)}, \tag 1$$
where $S_1,\dots,S_n$ are the left creation operators
on the full Fock space \linebreak $\F^2=\F^2(\H_n)$ (we refer to Section 1 for notation
and background material).
%
%

As  in [Po5], the noncommutative disc algebra $\Al_n$  is the
norm closed subalgebra  in $B(\F^2)$ generated by $S_1, \dots, S_n$ and the identity, and the Hardy (noncommutative analytic Toeplitz) algebra
  $F^\infty$ is the WOT-closed algebra generated by $\Al_n$ in
$B(\F^2)$.

It was proved in [Po8] that if $T=[T_1,\dots,T_n]$ is a contraction, then the map
 $$\Phi_T:C^*(S_1,\dots,S_n)\to B(\H);\quad
\Phi_T(S_{i_1}\cdots S_{i_k} S_{j_1}^*\cdots S_{j_p}^*)=
T_{i_1}\cdots T_{i_k} T_{j_1}^*\cdots T_{j_p}^*,$$
$1\le i_1,\dots,i_k,j_1,\dots,j_p\le n$,
is a completely contractive linear map, and $\Phi_T|_{\Al_n}$ is a homomorphism.
An elementary proof of this as well as an extension to a more general setting was obtained in [Po9], by the second author, using noncommutative Poisson transforms
on $C^*$-algebras generated by isometries (we refer to Section 3 for a sketch of the proof).

Let $J$ be a closed, $2$-sided ideal of $\Al_n$ with 
$J\subset \text{\rm Ker\,}\Phi_T$ and let $\N_J$ be the orthogonal of the image of $J$ in $\F^2$.  For each $i=1,\dots,n$, let 
$B_i:=P_{\N_J} S_i|_{\N_J}$.
Using noncommutative Poisson transforms [Po9], we will prove in Section 3 
that, for a large class of row contractions 
$T=[T_1,\dots,T_n]$ (including $C_0$-contractions), there is a unital, completely contractive, linear map 
$\Phi:C^*(B_1,\dots,B_n)\to B(\H)$ such that
$$
\Phi(B_{i_1}\cdots B_{i_k} B_{j_1}^*\cdots B_{j_p}^*)=
T_{i_1}\cdots T_{i_k} T_{j_1}^*\cdots T_{j_p}^*,$$
$1\le i_1,\dots,i_k,j_1,\dots,j_p\le n$.

The noncommutative dilation theory for $n$-tuples of 
operators [Fr], [Bu], [Po1], [Po2] was used in [Po6] to obtain
an $F^\infty$-functional calculus associated 
to any completely non-coisometric contraction (in short c.n.c.) $T=[T_1,\dots,T_n]$. More precisely, it was shown that the map 
$
\Psi_T: F^\infty\to B(\H)
$
defined by
$$
\Psi_T(f)=f(T_1,\dots,T_n):=\text{\rm SOT-}\lim_{r\to 1} f(rT_1,\dots,rT_n)
$$
is a
WOT-continuous and completely contractive homomorphism.
We will show that if $J$ is a WOT-closed, 2-sided ideal of 
$F^\infty$ with $J\subset \text{\rm Ker\,}\Psi_T$, then the map
$$
p(B_1,\dots,B_n)\mapsto p(T_1,\dots,T_n)
$$
can be extended to a WOT-continuous, completely contractive homomorphism from
$\W(B_1,\dots,B_n)$, the WOT-closed algebra generated by the compressions \linebreak $B_i:=P_{\N_J} S_i|_{\N_J}, \  i=1,\dots,n,$ to $B(\H)$.

Let us recall that $F^\infty$ has $\Bbb A_1(1)$ property (see [DP1]), therefore the $w^*$ and WOT topologies coincide on $F^\infty$.
An important step in proving the above-mentioned results is an extention of 
Sarason's result [S] to $F^\infty$. More precisely, we will show that if
 $J$ is a WOT-closed, 2-sided ideal of $F^\infty$, then
the map $$\Psi:F^\infty/J\to B(\N_J),\quad \Psi(f+J)=
P_{\N_J} f {}\vert_{\N_J}$$ is a $w^*$-continuous, completely isometric
representation.  In particular, for every $f\in F^\infty$,
$$\text{\rm dist}(f,J)=\|P_{\N_J} f {}\vert {\N_J}\|.$$

We present in this paper two proofs of this result: one is based on the noncommutative commutant lifting theorem [Po3] (see [SzF] and [FFr] for the classical case) and the characterization of the commutant of $S_1,\dots,S_n$ from [Po7], and the other is  based on noncommutative Poisson transforms [Po9] and representations of quotient algebras.
This is a key result which   
leads to the  noncommutative  interpolation theorems of 
Caratheodory (obtained previously
in [Po7]) and of Nevanlinna-Pick in the 
noncommutative analytic Toeplitz algebra $F^\infty$. 
Let us mentioned just one consequence of our results to the interpolation by bounded analytic functions in the unit ball of $\Bbb C^n$.
We will show that if
 $\lambda_1,\dots \lambda_k$ are $k$ distinct points in $\Bbb B_n$, the open unit ball
of $\Bbb C^n$,   $W_1,\dots,W_k\in B(\K)$ ($\K$ is a Hilbert space), and
the operator  matrix
$$
\left[{I_\K-W_j W_i^*\over 1-\langle \lambda_j,\lambda_i\rangle}\right]_{i,j=1,\dots,k}
$$
is positive definite, then there is an operator-valued analytic function
 \linebreak
$F:\Bbb B_n\to B(\K)$ such that 
$$\sup_{\zeta\in \Bbb B_n} \|F(\zeta)\|\le 1 \ \text{ and }
 \ F(\l_i)=W_j $$
for any  $j=1,\dots,k.$

In fact, we obtain more general results of interpolation by elements in $F^\infty$ (resp. $\Al_n$) with consequences to the interpolation by bounded analytic functions in the unit ball of $\Bbb C^n$.

\bigskip

\leftline{\bf 1. Notation and preliminary results }

Unless explicitly stated, $n$ stands for a cardinal number between
$1\leq n\leq\aleph_0$.  
Let $\H_n$ be an $n$-dimensional Hilbert space with orthonormal basis
$e_1,e_2,\dots,e_n$.   We consider the Full Fock space [E] of $\H_n$  
$$\F^2=\F^2(\H_n)=\bigoplus_{k\geq 0}\H_n^{\otimes k},$$ 
where $\H_n^{\otimes 0}=\C 1$ and $\H_n^{\otimes k}$ is the (Hilbert)
tensor product of $k$ copies of $\H_n$. 
We shall denote by $\P$ the set of all $p\in\F^2(\H_n)$ of the form
$$p=a_0  + \sum_{1\leq i_1,\dots, i_k\leq n\atop 1\leq k\leq m}
a_{i_1\cdots i_k}e_{i_1}\otimes e_{i_2}\otimes\cdots\otimes e_{i_k},\quad
m\in\Na,$$
where $a_0, a_{i_1\cdots i_k}\in\C$.  The set $\P$ may be viewed as the
algebra of the polynomials in $n$ noncommuting indeterminates, with
$p\otimes q$, $p,q\in\P$, as multiplication.  For any
  bounded operators $T_1,\dots,T_n$ on a Hilbert space $\H$,
define $$p(T_1,\dots,T_n):=a_0I_\H+\sum a_{i_1\cdots i_k}T_{i_1} T_{i_2}
\cdots T_{i_k}.$$

Let $\Fr_n^+$ be the unital free semigroup on $n$ generators $g_1,\dots,g_n$ and the identity $e$.
For each $\alpha\in\Fr_n^+$, define
$$e_\alpha := \cases e_{i_1}\otimes e_{i_2}\otimes \cdots \otimes e_{i_k},
&\text{if }\alpha=g_{i_1}g_{i_2}\cdots g_{i_k},\\
1, &\text{if }\alpha=e.\\
\endcases$$
It is easy to see that 
$\{e_\alpha:\alpha\in\Fr_n^+\}$ is an orthonormal basis of $\F^2$.
We also use $\Fr_n^+$ to denote arbitrary products of 
operators.  If $T_1,\dots,T_n\in B(\H)$, define 
$$T_\alpha := \cases T_{i_1}T_{i_2}\cdots T_{i_k},
&\text{if }\alpha=g_{i_1}g_{i_2}\cdots g_{i_k},\\
I_\H, &\text{if }\alpha=e.\\
\endcases$$
The length of $\alpha\in\Fr_n^+$ is defined by
$|\alpha|=k$, if $\alpha=g_{i_1}g_{i_2}\cdots g_{i_k}$, and
$|\alpha|=0$, if $\alpha=e$.
For each $i=1,\dots, n$, the left creation operator 
$$S_i:\F^2\to\F^2\quad \text{ is defined by }\quad S_i\psi=e_i\otimes\psi,\ \psi\in \F^2.$$
It is easy to see that $S_1,\dots,S_n$ are isometries with orthogonal
ranges. As in [Po5], 
$\Al_n$ is the norm closure of the algebra generated by 
$S_1,\dots,S_n$ and $I_{\F^2}$, and
$F^\infty$ is the weak operator topology closure of $\Al_n$.  
Alternatively, we let $\F^\infty$ be the set of those $\varphi\in\F^2$
such that
$$\|\varphi\|_\infty:=\sup\{\|\varphi\otimes p\|_2 : 
p\in\P, \|p\|_2\leq1\} <\infty.$$
For $\varphi\in\F^\infty$, define $\varphi(S_1,\dots,S_n):\F^2\to\F^2$
by $\varphi(S_1,\dots,S_n)\psi=\varphi\otimes\psi$.  
The norm $\|\varphi\|_\infty$ coincides with the operator norm of
$\varphi(S_1,\dots,S_n)$.  It will be useful later to view  
$\varphi\in\F^\infty$ as being an element in $F^\infty$ and conversely.
With this identification, $\Al_n$ is the closure of $\P$ in the
$\|\cdot\|_\infty$-norm.

We need to recall from [Po7] the characterization of the commutant of $S_1,\dots,S_n$.
Define the flipping operator $U:\F^2\to\F^2$ by
$$U (e_{i_1}\otimes e_{i_2}\otimes\cdots\otimes e_{i_k})=
 e_{i_k}\otimes\cdots\otimes e_{i_2}\otimes e_{i_1},$$ and let  $\tilde{\varphi}:=U\varphi$.  
It is easy to see that $U$ is a unitary operator, which satisfies
$U(\varphi\otimes \psi)=\tilde{\psi}\otimes\tilde{\varphi}$, and $U^2=I$.
An operator $A\in B(\F^2)$ commutes with $S_1,\dots,S_n$  if and only if
there exists $\phi\in \F^\infty$ such that
$
Ah=h\otimes\tilde{\phi},\quad h\in \F^2.
$
Notice that $A=U^*\phi(S_1,\dots,S_n)U$.

In [Po2], the second author defined $\varphi\in\F^\infty$ to be inner
if $\varphi(S_1,\dots,S_n)$ is an isometry, and outer if
$\varphi(S_1,\dots,S_n)$ has dense range. 
 A complete description of the invariant subspace structure
of $\F^\infty$ was obtained in [Po2] (even in a more general setting), 
using a noncommutative version of the Wold decomposition (see [Po1]).
A family of inner operators $\{\varphi_i:i\in I\}$ is called
orthogonal, if whenever $i\not=j$, $\F^2\otimes\tilde{\varphi}_i$ 
is orthogonal to $\F^2\otimes\tilde{\varphi}_j$; or equivalently,
$\varphi_i\otimes\F^2$ is orthogonal to $\varphi_j\otimes\F^2$.  
It follows from
[Po2; Theorem 2.2] that a subspace $\M$ of $\F^2$ is invariant
under $S_1,\dots,S_n$ if and only if $\M=\bigoplus_{i\in I}\F^2\otimes
\tilde{\varphi}_i$, for some family of orthogonal inner operators.

 The second author obtained in [Po4] an inner-outer factorization  
which implies that 
any $\eta\in\F^\infty$ can be factored as $\eta=\varphi\otimes\psi$,
where $\varphi$ is inner and $\psi$ is outer.  The same factorization
result was proved for elements of $\F^2$ in [APo], where 
$\eta\in\F^2$ was said to be outer if there exists 
a sequence of polynomials $p_n\in\P$ such that 
$\psi\otimes p_n\to 1$ in the norm of $\F^2$ (this last result was
also obtained  recently in [DP1]).
 Let us mention
that we proved in [APo] that the noncommutative analytic Toeplitz algebra in $n$ noncommuting variables
$\F^\infty$ is reflexive.
Recently, Davidson and Pitts [DP1] proved that this 
algebra is hyper-reflexive.
They also studied in [DP2] the algebraic structure
of $\F^\infty$  (in their notation $\L_n$).

Now let us recall some general facts about duality in Banach spaces.  
Let $X$ be a Banach space with predual $X_*$ and dual $X^*$, and let
$\S\subset X$.  The preannihilator of $\S$ in $X_*$ is the set
${}^\perp\S=\{f\in X_*: \langle f,x\rangle=0$ for all $x\in\S\}$, and
the annihilator of $\S$ in $X^*$ is the set
$\S^\perp=\{f\in X^*: \langle x,f\rangle=0$ for all $x\in\S\}$.
If $\S$ is $w^*$-closed, it is well known that 
$({}^\perp S)^*=X/S$, $({}^\perp\S)^\perp=\S$, $(X_*/{}^\perp\S)^*=\S$, and
$({}^\perp\S)^{**}=\S^\perp$ (see [Ru]).

The predual of $B(\H)$ is the space of trace class operators $c_1(\H)$,
under the trace duality.  That is, if $T\in c_1(\H)$ and $A\in B(\H)$, then
$\langle T,A\rangle :=tr(TA)$.  
A $w^*$-closed subspace $\S$ of $B(\H)$ has property $\A_1$ if for every 
$w^*$-continuous linear map $f:\S\to\C$, there exists $h, k\in\H$
such that for all $T\in\S$, $f(T)=\langle Th,k\rangle$.  Moreover, if
for every $\epsilon>0$, $h$ and $k$ can be chosen so that
$\|h\|\|h\|\leq (1+\epsilon)\|f\|$, $\S$ has property
$\A_1(1)$ (we refer to [BFP] for more information).
If $n\geq2$, Davidson and Pitts [DP1]
proved that $F^\infty\subset B(\F^2)$ 
has property $\A_1(1)$, and Bercovici [B] proved that
$M_k(F^\infty)$ has property $\A_1(1)$ for each $k\geq1$ 
(i.e., $F^\infty$ has property $\A_{\aleph_0}(1)$).

We refer to [Ar1], [P], and [Pi] for results on completely
bounded maps and operator spaces.

\medbreak

Let $J$ be a WOT-closed, 2-sided ideal of $F^\infty$ and define
$$\eqalign{\M_J 
=& \overline{\{\varphi\otimes\psi:\varphi\in J, \psi\in\F^2\}}^{\|\cdot\|_2},\text{ and }\cr
            \N_J=&\F^2\ominus\M_J.\cr}$$
\medbreak

Recall from [S] that a subspace $\N\subset\H$ is semi-invariant under
a semigroup of operators $\Sigma\subset B(\H)$ 
if for every $T_1,T_2\in \Sigma$,
$P_\N T_1  P_\N T_2  P_\N=
P_\N T_1 T_2 P_\N.$  It is
well known that if $\N_1,\N_2$ are invariant subspaces under $\Sigma$ and
$\N_2\subset\N_1$ then $\N_1\ominus\N_2$ is semi-invariant under 
$\Sigma$.

\proclaim{Lemma 1.1} If $J$ is a  \text{\rm WOT}-closed, $2$-sided ideal of $F^\infty$, then the subspaces $\N_J$ and $U\N_J$ are invariant under $S_i^*,\ i=1,2,\dots,n$.
\endproclaim

\demo{Proof}
Since $J$ is a left ideal, $\M_J$ is invariant under $S_1,\dots,S_n$ and,
hence, $\N_J$ is invariant under $S_1^*,\dots,S_n^*$.  Moreover, since $J$ is 
a right ideal, the set \linebreak $\{\varphi: \varphi(S_1,\dots,S_n)\in J\}$ is dense
in $\M_J$. Similarly, one can prove that  $U\N_J$ is invariant to each
$S_i^*,\ i=1,2,\dots,n$. $\blacksquare$
\enddemo

\proclaim{Proposition 1.2} Let $J$ be a $w^*$-closed, $2$-sided ideal 
of $F^\infty$ and $g\in {}^\perp J\subset (F^\infty)_*$.
For each $\epsilon>0$, there exist $\psi_1,\psi_2\in\N_J$ satisfying
$\|\psi_1\|_2 \|\psi_2\|_2\leq (1+\epsilon)\|g\|$ such that for every 
$\eta\in F^\infty$,  $g(\eta)=\langle\eta\otimes\psi_1,\psi_2\rangle$. Conversely, if
$\psi_1,\psi_2\in\N_J$ and $g(\eta):=\langle\eta\otimes\psi_1,\psi_2\rangle$ for any  $\eta\in F^\infty$, then  $g\in {}^\perp J$.
\endproclaim

\demo{Proof}  Let $g\in {}^\perp J\subset (F^\infty)_*$ and let $\epsilon>0$.
Since $F^\infty$ has the $\A_1(1)$ property, find
$\varphi_1,\varphi_2\in\F^2$ satisfying
$\|\varphi_1\|_2\|\varphi_2\|_2\leq (1+\epsilon)\|g\|$ such that 
for every $\xi\in F^\infty$,
$\langle g, \xi\rangle=\langle\xi\otimes\varphi_1,\varphi_2\rangle$.  Factor $\varphi_1=\eta_1\otimes\eta_2$,
where $\tilde{\eta}_1$ is outer and $\tilde{\eta}_2$ is inner.  Hence,
there exist a sequence of polynomials $p_n\in\P$ such that 
$p_n\otimes\eta_1\to 1$ in the norm of $\F^2$, and $\F^2\otimes\eta_2$
is a closed subspace of $\F^2$.  Let $P_{\eta_2}$ be the orthogonal
projection onto $\F^2\otimes\eta_2$, and write $P_{\eta_2}\varphi_2=
\psi_2\otimes\eta_2$ for some $\psi_2\in\F^2$.  
Then for each $\xi\in\F^\infty$,
$$\eqalign{g(\xi)
&=\langle\xi\otimes\varphi_1,\varphi_2\rangle
	=\langle \xi\otimes\eta_1\otimes\eta_2,\varphi_2\rangle 
	 =\langle P_{\eta_2} (\xi\otimes\eta_1\otimes\eta_2),\varphi_2\rangle\cr
&=\langle \xi\otimes\eta_1\otimes\eta_2,P_{\eta_2} (\varphi_2) \rangle
	 =\langle \xi\otimes\eta_1\otimes\eta_2,\psi_2\otimes\eta_2\rangle
		=\langle \xi\otimes\eta_1,\psi_2\rangle.\cr}$$
The last equality follows because the operator on $\F^2$
that multiplies from the right by
$\eta_2$ is an isometry (this is the equivalent to $\tilde{\eta}_2$ inner).

We will show that $\psi_2\in\N_J$.  
Let $\xi\in J$.  Since $J$ is a right ideal,
$\xi\otimes p_n\in J$ for each $n$.  Hence,
$0=g(\xi\otimes p_n)=\langle \xi\otimes p_n\otimes \eta_1,\psi_2\rangle
\to \langle\xi,\psi_2\rangle.$ Therefore, $\langle\xi, \psi_2\rangle=0$.
Since $\{\varphi\in\F^2 : \varphi(S_1,\dots, S_n)\in J\}$ is dense in
$\M_J$ we conclude that $\psi_2\in\N_J$.

Recall that $\N_J$ is invariant under $S_\alpha^*$ and 
let $\psi_1=P_{\N_J}(\eta_1)$.  For each $\alpha\in\Fr_n^+$,
  $$
\align
g(e_\alpha)&=\langle e_\alpha\otimes\eta_1,\psi_2\rangle=
\langle\eta_1,S_\alpha^*(\psi_2)\rangle\\
&=\langle\eta_1,P_{\N_J}S_\alpha^*(\psi_2)\rangle=
\langle P_{\N_J}(\eta_1),S_\alpha^*(\psi_2)\rangle\\
&=\langle\psi_1,S_\alpha^*(\psi_2)\rangle
=\langle e_\alpha\otimes\psi_1,\psi_2\rangle.
\endalign
$$
The converse is straightforward.
  This completes the proof.
$\blacksquare$
\enddemo
\medbreak

As a consequence of Proposition 1.2, we obtain the following.

\proclaim{Proposition 1.3} For every $\varphi \in F^\infty, \text{\rm dist}
(\varphi,J)=
\|P_{\N_J}\varphi(S_1,\dots,S_n)\vert_{\N_J}\|$.
Consequently, the map $\Phi:F^\infty/J\to B(\N_J)$ defined by
$\Phi(\varphi+J)= P_{\N_J}\varphi(S_1,\dots,S_n)\vert_{\N_J}$
is an isometric homomorphism.
\endproclaim

\demo{Proof}  If $\psi\in J$ it is clear that 
$P_{\N_J}\psi(S_1,\dots,S_n)\vert_{N_J}=0$.  Hence, for every
$\varphi\in F^\infty$,
$\|P_{\N_J}\varphi(S_1,\dots,S_n)\vert_{\N_J}\|\leq\text{dist}(\varphi,J)$.

Suppose now that $\varphi\not\in J$.  
Since $({}^\perp J)^*=F^\infty/J$, for 
every $\epsilon>0$ there exists $f\in {}^\perp J$, $\|f\|<1$ such that
$|f(\varphi)|>\text{dist}(\varphi,J)-\epsilon$.  By Proposition 1.2, there
exist $\xi_1, \xi_2\in\N_J$ such that $\|\xi_1\|_2 \|\xi_2\|_2\leq 1$
and $f(\varphi)=\langle \varphi(S_1,\dots,S_n)\xi_1,\xi_2\rangle$.
Hence, $\|P_{\N_J}\varphi(S_1,\dots,S_n)\vert_{\N_J}\|\geq
\text{dist}(\varphi,J)-\epsilon$.  Since $\epsilon>0$ is
arbitrary, we finish the proof.\ 
$\blacksquare$
\enddemo
 
It should be noted that Proposition 1.3 is all one really needs to 
obtain the scalar version of Caratheodory or
Nevanlinna-Pick interpolation in $F^\infty$.

\medbreak

We know from [Po5], [Po6] that the set $\P$ of all polynomials in $S_1,\dots, S_n$ is WOT-dense in $\F^\infty$. Indeed, if
       $f=\sum_{\alpha\in \Bbb F_n^+} a_\alpha e_\alpha$ is
 in $\F^\infty$ and 
  $f_r:=\sum_{\alpha\in \Bbb F_n^+}r^{|\alpha|} a_\alpha e_\alpha$
  for any $0<r<1$ then SOT-$\lim_{r\to 1}f_r=f$ and $\|f_r\|_\infty
  \le\|f\|_\infty$ (see [Po6]). On the other hand, $f_r\in\Al_n$.
Indeed,
  since $\|\sum_{|\alpha|=k}a_\alpha S_\alpha \|=(\sum_{|\alpha|=k}|a_\alpha
  |^2)^{1/2}$, we have
  $$
  \align
  \sum_{k=0}^\infty r^k\|\sum_{|\alpha|=k}a_\alpha S_\alpha \|
  &=  \sum_{k=0}^\infty r^k(\sum_{|\alpha|=k}|a_\alpha
  |^2)^{1/2}\\
  &=(\sum_{k=0}^\infty r^k)\|\phi\|_2.
  \endalign
  $$
   Therefore, $\sum_{\alpha\in \Bbb F_n^+}r^{|\alpha|} a_\alpha S_\alpha$
   converges in norm, so that $f_r\in\Al_n$. Taking into account that $\P$ is norm dense in $\Al_n$, the result follows.

  We will use the same notation as above if $J$ is a closed, $2$-sided ideal in $\Al_n$.
\proclaim{Lemma 1.4}                      
Let $J\subset \Al_n$ be a  $2$-sided ideal of $\Al_n$ and let $J_w$ be the \text{\rm WOT}-closed, $2$-sided ideal  generated by $J$ in $F^\infty$.
Then $J_w={\bar J}^{\text WOT}$ and   $\N_{J_w}=\N_J$.
\endproclaim
\demo{Proof}
  We need to show that ${\bar J}^{\text{\rm WOT}}$ is a 2-sided ideal in $F^\infty$.
    Consider \linebreak $\psi,\phi\in F^\infty,\ f\in {\bar J}^{\text WOT}$, and let
    $\{g_i\}\subset J$ be a net WOT-convergent to $f$.

  Since $J$ is a 2-sided ideal of $\Al_n$, $\phi_r g_i\psi_{r'}\in J$ \ for any
 $r,r'\in (0,1)$ and any $i$.
 Using the remarks preceding this lemma, it is easy to see, by taking appropriate limits, that
 $\phi f \psi\in{\bar J}^{\text WOT}$.
 Now let us show that $\N_{J_w}=\N_J$. Since $J\subset J_w$, it is clear that
 $\N_J\supset\N_{J_w}$.
 Let $f\in{\bar J}^{\text WOT},\ \psi\in\F^2$, and choose $\{g_i\}\subset J$ such that
 WOT-$\lim_i g_i=f$.
 If $x\in\N_J$ we have
 $$
 \langle x,f\otimes\psi\rangle=\lim_i
 \langle x, g_i(S_1,\dots,S_n)\psi\rangle=0.
 $$
Therefore, $x\in\N_{J_w}$, which proves that $\N_J\subset\N_{J_w}$.
This completes the proof.
$\blacksquare$
\enddemo
%
%
\bigskip

 \leftline{\bf 2. Non-commutative   interpolation in $F^\infty$}

    Let $\H, \K$ be Hilbert spaces and $I$ be a set of indices with
    $\dim\K=\text{\rm card\,} I=\ga$.
    Denote $\oplus_\ga\H:=\oplus_{i\in I}\H_i$ where $\H_i:=\H$,
    and notice that, under the canonical identification
    $\oplus_\ga\H=\H\otimes\K$ (Hilbert tensor product), each operator 
    $X\in B(\H\otimes\K)$ can be seen as a {\it matrix } of operators in $B(\H)$,
    i.e., $X=[X_{\alpha\beta}]_{\alpha,\beta\in I}$ with \linebreak
    $X_{\alpha\beta}\in B(\H)$.
    For any $\U\subset B(\H)$, we denote
    $$
   M_\ga(\U)=\{[u_{\alpha\beta}]\in B(\oplus_\ga\H)\ : \ u_{\alpha\beta}\in\U; {\alpha,\beta}\in I\}.
    $$
    It is clear now that
    $
    \{u\otimes I_\K;\ u\in\U\}'=M_\ga(\U')
    $
    (where $'$ stands for commutant).

Let us recall from [Po1] that $T=[T_1,\dots,T_n]$ is  called $C_0$-contraction if $T$ is a contraction and 
$$\text{\rm SOT-}\lim_{k\to\infty}\sum_{\alpha\in \Bbb F_n^+,|\alpha|=k}T_\alpha T_\alpha^*=0.\leqno{(2.1)}$$
For example, if $\sum_{i=1}^n T_iT_i^*\le \rho I_\H$ for some $0<\rho<1$, then $[T_1,\dots,T_n]$ is a $C_0$-contraction.

The following result is an extension of Sarason's result [S] and a consequence
of the noncommutative commutant lifting theorem [Po3] and the characterization of the commutant of $\{S_1,\dots,S_n\}$ from [Po7].
    \proclaim{Theorem 2.1}
    Let $\K$ be a Hilbert space with $\dim \K=\ga$ and let $\N \subset\F^2$ be an invariant subspace
    for $S^*_1,\dots,S_n^*$.
    If $T\in B(\N\otimes\K)$ commutes with each \linebreak
    $X_i:= P_\N S_i|_\N\otimes I_\K, \ i=1,2,\dots,n$, then there is
    $\Phi(S_1,\dots,S_n):=[\phi_{\alpha,\beta}(S_1,\dots,S_n)]$ in $M_\ga(F^\infty)$
    such that $\|\Phi(S_1,\dots,S_n)\|=\|T\|$ and
    $$
    P_{\N\otimes\K}[U^*\phi_{\alpha\beta}(S_1,\dots,S_n)U]=T P_{\N\otimes\K},
    $$
    where $P_{\N\otimes\K}$ is the orthogonal projection of $\F^2\otimes\K$ onto
    $\N\otimes\K$, and $U$ is the flipping operator on $\F^2$.
    \endproclaim
    \demo{Proof}
    Since $S_i^*|_\N=B_i^*, i=1,2,\dots,n$, it is clear that 
    $[P_\N S_1|_\N,\dots, P_\N S_n|_\N]$ is a $C_0$-contraction and, according to
    [Po1], its minimal isometric dilation is $[S_1,\dots,S_n]$.
    Therefore, the minimal isometric dilation of $[X_1,\dots,X_n]$ is
    $[S_1\otimes I_\K,\dots,S_n\otimes I_\K]$.
    According to  the noncommutative commutant lifting theorem [Po3], there is \linebreak $A\in\{S_i\otimes I_\K; i=1,2,\dots,n\}'$
    such that $\|A\|=\|T\|$ and $A^*|_{\N\otimes\K}=T^*$.
    Therefore, there exists $a_{\alpha,\beta}\in\{S_1,\dots,S_n\}'$
    such that $A=[a_{\alpha\beta}]\in M_\ga(\{S_1,\dots,S_n\}')\subset
    B(\oplus_\ga \F^2)$.
     Using the characterization of the commutant of $\{S_1,\dots, S_n\}$ from [Po7], we find
     $\phi_{\alpha\beta}\in\F^\infty$
    such that $a_{\alpha\beta}= U^*\phi_{\alpha\beta}(S_1,\dots,S_n)U$,
    where $U$ is  \linebreak  the flipping operator on $\F^2$.
    Therefore $A=[U^*\phi_{\alpha\beta}(S_1,\dots,S_n)U]$ and \linebreak
    $T=[P_\N U^*\phi_{\alpha\beta}(S_1,\dots,S_n)U|_\N]$.
        $\blacksquare$
\enddemo

Notice that if $n=1$ we find again Sarason's result [S].

\proclaim{Lemma 2.2}
    Let $T_i\in B(\H)$ be such that $T:=[T_1,\dots,T_n]$ is a $C_0$-contraction and let
    $f_{\alpha\beta}(S_1,\dots,S_n)\in F^\infty, \alpha,\beta\in I$, be such
    that
    $[f_{\alpha\beta}(S_1,\dots,S_n)]_{\alpha,\beta\in I}$ is in $B(\oplus_\ga\F^2)
    \quad(\ga=\text{\rm card\,} I)$.
    Then $[f_{\alpha\beta}(T_1,\dots,T_n)]\in B(\oplus_\ga\H)$
    and
    $$\|[f_{\alpha\beta}(T_1,\dots,T_n)]\|\leq\|[f_{\alpha\beta}(S_1,\dots,S_n)]\|.$$
    \endproclaim
    \demo{Proof}
    According to [Po1], the minimal isometric dilation of
    $T=[T_1,\dots,T_n]$ is
    $[S_1\otimes I_\L,\dots, S_n\otimes I_\L]$ for some Hilbert space $\L$.
   According to
    Theorem 3.6 from [Po5], for any
    $f_{\alpha\beta}\in\F^\infty$, 
    $
    f_{\alpha\beta}(T_1,\dots,T_n)=P_\H(f_{\alpha\beta}(S_1,\dots,S_n)
    \otimes I_\L)|_\H.
    $
    The rest of the proof is straightforward.
        $\blacksquare$
\enddemo

If $\zeta=(\zeta_1,\dots,\zeta_n)\in \Bbb C^n$ is such that $|\zeta|:=(|\zeta_1|^2+\cdots+\zeta_n|^2)^{1/2}<1$, and 
$f(S_1,\dots, S_n)\in F^\infty$, then, according to the $F^\infty$-functional calculus [Po6], we infer that 
$|f(\zeta_1,\dots,\zeta_n)|\le \|f(S_1,\dots, S_n)\|$.
Therefore, $f(\zeta_1,\dots,\zeta_n)$ is an analytic function in
$\Bbb B_n$.
Moreover, we deduce the following.
    \proclaim{Corollary 2.3}
    If $f_{\alpha\beta}\in F^\infty, \ \alpha,\beta\in I$,  $\text{\rm card\,}I=\ga$, and
    $[f_{\alpha\beta}(S_1,\dots,S_n)]_{\alpha,\beta\in I}$ is in
    $B(\oplus_\ga\F^2)$,
    then
    $\Phi(\zeta):=[f_{\alpha\beta}(\zeta)]$ is an operator-valued
    analytic function in $\Bbb B_n$.
    Moreover, $\sup_{\zeta\in \Bbb B_n}\|\Phi(\zeta)\|\leq
    \|[f_{\alpha\beta}(S_1,\dots,S_n)]\|$.
    \endproclaim

A consequence of Theorem 2.1 is the following extension of the 
Nevanlinna-Pick problem to the noncommutative Toeplitz algebra $F^\infty$.

\proclaim{Theorem 2.4}
 Let $\la_1,\dots,\la_k$ be $k$ district points in $\Bbb B_n$ and let
    $W_1,\dots,W_k$ be in $B(\K)$, where $\K$ is a Hilbert space with $\dim\K=\ga$.
    Then there exists
    $\Phi(S_1,\dots,S_n):=[\phi_{\alpha,\beta}(S_1,\dots,S_n)]$ in
    $M_\ga(F^\infty)$, such that $\|\Phi(S_1,\dots,S_n)\|\leq 1$ and
    $\Phi(\la_j)=W_j$,  $j=1,2,\dots,k$, if and only if the operator matrix
    $$
    \[{I_\K-W_j W_i^*\over 1-\langle\la_j,\la_i\rangle}\]_{i,j=1,2,\dots,k} \tag 2.2
    $$
    is positive definite.
    \endproclaim

\demo{Proof} For each $i=1,\dots,k$, let $\l_i:=(\l_{i1},\dots,\l_{in})\in \Bbb C^n$ and, for $\alpha=g_{j_1}g_{j_2}\dots g_{j_m}$ in $\Bbb F_n^+$,
let $\l_{i\alpha}:=\l_{ij_1}\l_{ij_2}\dots\l_{ij_m}$ and $\l_e=1$.
Define 
$$
z_{\l_i}:=\sum_{\alpha\in \Bbb F_n^+} \overline {\l}_{i\alpha} e_\alpha,\quad 
i=1,2,\dots,n,
$$
and notice that, for any $\phi=\sum_{\alpha\in \Bbb F_n^+}a_\alpha e_\alpha$ in $\F^2$,\  $\langle\phi, z_{\l_i}\rangle=\phi(\l_i)$. 

If $\phi\in\F^\infty$ then $\langle\phi, z_{\l_i}\rangle=\langle 1,\phi(S_1,\dots , S_n)^* z_{\l_i}\rangle$, where $S_1,\dots , S_n$ are the left
creation operators on the full Fock space $\F^2$.
It is clear that $S_i^* z_{\l_j}=\overline{\l}_{ji} z_{\l_j}$
for any $i=1,\dots,n;\ j=1,\dots,k$.
Denote 
$$
\N:=\text{span} \{ z_{\l_j}:\ j=1,\dots,k\}
$$
  and define $X_i\in B(\N\otimes\K)$ by $X_i=P_\N S_i|_\N\otimes I_\K$.
    Since $z_{\la_1},\dots, z_{\la_k}$ are linearly independent, we can define
    $T\in  B(\N\otimes\K)$ by setting
    $$
    T^*(z_{\la_j}\otimes h)=z_{\la_j}\otimes W_j^* h
\tag 2.3    $$
    for any $h\in \H,\  j=1,\dots, k$.
    Notice that for each $i=1,\dots, k $, $TX_i=X_iT$. Indeed,
     $$
     \align
     X_i^*T^*(z_{\la_j}\otimes h)&=X_i^*(z_{\la_j}\otimes W_j^* h)=
     S_i^*z_{\la_j}\otimes W_j^* h\\
     &=\overline{\la}_{ji}z_{\la_j}\otimes W_j^* h
     \endalign
     $$
     and
     $$  T^*X_i^*(z_{\la_j}\otimes h)=T^*(\overline{\la}_{ji}z_{\la_j}\otimes h)=
     \overline{\la}_{ji}z_{\la_j}\otimes W_j^* h.
     $$
      Since $\N$ is invariant under $S_i^*, \ i=1,\dots,n$, according to Theorem
      2.1, there exists $\[\phi_{\alpha,\beta}(S_1,\dots, S_n)\]\in
       M_\ga( F^\infty)$ such that $P_{\N\otimes \K}\[U^*
       \phi_{\alpha,\beta}(S_1,\dots, S_n)U\]=TP_{\N\otimes \K},$
       and
       $$\|\[\phi_{\alpha,\beta}(S_1,\dots, S_n)\]\|=\|T\|.\tag 2.4$$
    Let us show that $\Phi(S_1,\dots, S_n):=\[\phi_{\alpha,\beta}(S_1,
    \dots, S_n)\]$ satisfies $\Phi(\la_j)=W_j$ for any $j=1,\dots,k$,
    if and only if
    $$  \[P_\N U^*\phi_{\alpha,\beta}(S_1,\dots, S_n)U|_\N\]=T.\tag 2.5
    $$
    To prove this, notice first that $U(z_{\la_j})=z_{\la_j}, \ j=1,\dots,k$,
    and
    $$
    \langle\phi_{\alpha,\beta}(S_1,\dots, S_n) z_{\la_j}, z_{\la_j}\rangle=
     \phi_{\alpha,\beta}(\la_j)  \langle     z_{\la_j}, z_{\la_j}\rangle.
    $$
    Due to these relations, (2.3), and (2.5), it is easy to see that, for any
    $j=1,\dots,k$ and $h,k\in \K$, we have
    $$
    \align
    \langle [U^*\phi_{\alpha,\beta}(S_1,\dots, S_n)U ]&
    ( z_{\la_j} \otimes h), z_{\la_j}\otimes k\rangle =
     \langle z_{\la_j}, z_{\la_j}\rangle
      \langle \Phi(\la_j)h, k\rangle\\
      &= \langle T(z_{\la_j} \otimes h), z_{\la_j}\otimes k\rangle=
     \langle z_{\la_j} \otimes h, z_{\la_j}\otimes W_j^* k\rangle\\
     &=\langle z_{\la_j}, z_{\la_j}\rangle
      \langle W^*_j h, k\rangle.
      \endalign
      $$

      Now, it is clear that  $\Phi(\la_j)=W_j$ for any $j=1,\dots,k$ if
       and only if (2.5) holds.
       On the other hand, (2.4) shows that $\|\Phi(S_1,\dots,S_n)\|\leq 1$
       if and only if $\|T\|\leq 1$.
       The later condition is equivalent to
       $$
       \langle g\otimes k, g\otimes k\rangle -\langle T^*( g\otimes k),
       T^*( g\otimes k)\rangle\ge 0
       $$
 for any $g=\sum_{i=1}^k \alpha_j z_{\l_j}$ in $\N$ and $k\in \K$.
This inequality is equivalent to 
$$
\sum_{i,j=1}^k \alpha_i\overline{\alpha}_j \langle z_{\l_i}, z_{\l_j}\rangle(I-W_j W^*_i)\ge 0, \tag 2.6
$$
for any $\alpha_i\in\Bbb C$ and $k\in \K$.
Since 
 $$
\align
\langle z_{\l_i}, z_{\l_j}\rangle &=z_{\l_i}(\l_j)=\sum_{\alpha\in \Bbb F_n^+}
\overline{\l}_{i\alpha} \l_{j\alpha}\\
&=1+\langle\l_j,\l_i\rangle+\langle\l_j,\l_i\rangle^2+\dots\\
&={1\over 1-\langle\l_j,\l_i\rangle},
\endalign
$$
inequality (2.6) holds  if and only if the matrix (2.2) is positive definite.
This completes the proof.
$\blacksquare$
\enddemo

\proclaim{Corollary 2.5} If $n=1$ we find again the Nevanlinna-Pick interpolation theorem {\rm (see [Pic], [N])}.
\endproclaim

Notice that the proof of Theorem 2.1 works also for arbitrary families 
$\{\la_j\}_{j\in J}$ of distinct elements in $\Bbb B_n$, the open  unit
 ball of $\Bbb C^n$.

\proclaim{Theorem 2.6} Let $\{\la_j\}_{j\in J}$ be distinct elements
 in $\Bbb B_n$ and let $\{W_j\}_{j\in J}\subset B(\K)$, where $\K$
    is a Hilbert space of dimension $\ga$.
Then there exists $\Phi(S_1,\dots, S_n)$ in $M_\ga( F^\infty)$
 such that $ \|\Phi(S_1,\dots, S_n)\|\le 1$ and  $\Phi(\la_j)=W_j$ for all
 $j\in J$ if and only if
$$
\sum_{i,j\in J}  
{I_\K-  W_j W^*_i\over 1-\langle \lambda_j,\lambda_i\rangle}\alpha_i\overline{\alpha}_j\ge 0
$$
for any $\{\alpha_j\}_{j\in J}$ such that $\{j:\alpha_j\neq 0\}$ is finite.
\endproclaim

Combining  Theorem 2.4 with Corollary 2.3, we obtain the following
sufficient condition for interpolation in the open unit ball of $\Bbb C^n$.

\proclaim{Corollary 2.7}
Let $\lambda_1,\dots \lambda_k$ be $k$ distinct points in $\Bbb B_n$ and
 let  $W_1,\dots, W_k$ be in $B(\K)$.
If the matrix
$$
\left[{1-  W_jW^*_i\over 1-\langle \lambda_j,\lambda_i\rangle}\right]_{i,j=1,\dots,k}\tag2.7
$$
is positive definite, then there is an operator-valued analytic function
\linebreak  $F: \Bbb B_n\to B(\K)$   such that
  $$\sup_{\zeta\in \Bbb B_n} \|F(\zeta)\|\le 1 \ \text{ and }
 \ F(\la_i)=W_j$$
 for any  $j=1,\dots,k.$
\endproclaim

Arveson [Ar2] showed that there are functions $F$ in $H^\infty(\Bbb B_n)$ for
which there are no $f\in F^\infty$ such that $f(\lambda)=F(\lambda)$
for each $\lambda\in{\Bbb B}_n$.
The next result characterizes those functions in $H^\infty(\Bbb B_n)$
which are the image of elements in the unit ball of $F^\infty$.

\proclaim{Theorem 2.8}
Let $F$ be a complex-valued function defined on $\Bbb B_n$, such that \linebreak $|F(\zeta)|<1$  for all $|\zeta|<1$. Then there is  $f\in F^\infty,\ \|f\|_\infty\le 1$ such that \linebreak $f(\zeta)=F(\zeta),\ \zeta\in\Bbb B_n$, if and only if for each $k\ge 1$ and each $k$-tuple of points \linebreak $\la_1,\dots,\la_k\in \Bbb B_n$, the matrix 
$$
\left[{1-F(\la_j)\overline{F(\la_i)} \over 1-\langle \lambda_j,\lambda_i\rangle}\right]_{i,j=1,\dots,k}\tag 2.8
$$
is positive definite. In particular, if $(2.8)$ holds, then $F$ is analytic on $\Bbb B_n$. 
\endproclaim
\demo{Proof}
 The necessity of (2.8) follows immediately from Theorem 2.4.
 Conversely, suppose that $F$ satisfies (2.8). Let $\{\la_j\}_{j=1}^\infty$
 be a countable dense set in $\Bbb B_n$. According to Theorem 2.4, for each $k$,
  there is $f_k\in F^\infty$ with $\|f_k\|_\infty\le 1$ and
  $$
  f_k(\la_j)=F(\la_j)\quad \text{ for any } j=1,\dots,k.\tag 2.9
  $$
  Since  $\{f_k\}_{k=1}^\infty$ is bounded and $F^\infty $ is a dual space,
   according to Alaoglu's theorem, there is a subsequence
   $\{f_{k_m}\}_{k=1}^\infty$  such that $f_{k_m}$ converges in the  $w^*$-
   topology to an element $f\in F^\infty,\ \|f\|_\infty\le 1$.
   Since $w^*$ and WOT topologies coincide on $F^\infty$ and the $F^\infty$-
   functional calculus for $C_0$-contractions is WOT-continuous, we infer
   that
   $$
   \lim_{m\to \infty} f_{k_m}(\la_{j1},\dots, \la_{jn})=f(\la_{j1},\dots,
   \la_{jn}),\quad \text{ where } \la_j=(\la_{j1},\dots, \la_{jn}).
    $$
    Using (2.9), we have
    $\lim_{m\to \infty} f_{k_m}(\la_{j1},\dots, \la_{jn})=F(\la_{j1},\dots,
    \la_{jn})$.
    Therefore, $f(\la_j)=F(\la_j)$ for any $j=1,2, \dots.$

    We claim that $f(\zeta)=F(\zeta)$ for any $\zeta\in \Bbb B_n$.
    Let $\la$ be an arbitrary point in $\Bbb B_n$. By repeating
     the preceding argument, there is $g\in F^\infty, \ \|g\|\le 1$ so that
      $g(\zeta)=F(\zeta)$ on the set
       $\{\la_j\}_{j=1}^\infty\cup \{\la\}$.
       Since the maps $\zeta\mapsto g(\zeta)$ and   $\zeta\mapsto f(\zeta)$
     are analytic in $\Bbb B_n$ and coincide on  $\{\la_j\}_{j=1}^\infty$, which
     is dense in $\Bbb B_n$, we infer that they coincide on $\Bbb B_n$.
     In particular, we obtain $f(\la)=F(\la)$. Since $\la$ was an arbitrary point in
    in $\Bbb B_n$, we deduce that $f$ and $F$ coincide on $\Bbb B_n$.
    In particular, $\zeta\mapsto F(\zeta)$ is a bounded analytic function
    in $\Bbb B_n$. This completes the proof.
$\blacksquare$
    \enddemo

Condition 
$(2.7)$ is necessary and sufficient for interpolation in $F^\infty$ but only
sufficient for interpolation in $H^\infty(\Bbb B_n)$. 
One can  use the classical Cauchy formula for $\Bbb B_n$ to obtain a necessary condition
for Nevanlinna-Pick interpolation in $H^\infty(\Bbb B_n)$.  Recall that
for every $f\in  H^\infty(\Bbb B_n)$ and $\lambda\in \Bbb B_n$,
$$f(\lambda)=\int_{\partial \Bbb B_n} 
{f(w)\over (1-\langle \lambda,w\rangle)^n}d\sigma(w)$$
where $\sigma$ is the rotation invariant probability measure on 
$\partial \Bbb B_n$.
Using this formula, and a standard argument (eg. like the one used in
Section 3 of [CW]) we can check that if there exists 
$f\in H^\infty(\Bbb B_n)$, $\|f\|_\infty\leq1$, such that
$f(\lambda_j)=w_j$ for $j=1,\dots,k$, then
$$\left[{1-w_j{\bar w}_i \over (1-\langle\lambda_j,\lambda_i\rangle)^n}
\right]_{i,j=1,\dots,k}\tag 2.10$$
is positive definite.
Now, one can easily check that the scalar version of
condition (2.7) is stronger than condition (2.10)
(see for example Lemma 4.1 of [CW]).

\smallbreak

    Let $J$ be a $w^*$-closed, 2-sided ideal of $F^\infty$.
    For any cardinal $\ga$, the algebra $ M_\ga(F^\infty)$ is $w^*$-closed in
    $B(\oplus_\ga \F^2)$ and $M_\ga(J)$ is a $w^*$-closed, 2-sided ideal of
    $M_\ga(F^\infty)$.
Recently, Bercovici [B] proved that if the commutant
of a $w^*$-closed subspace of $B(\H)$ contains two isometries with
orthogonal ranges, then the subspace has property $X_{0,1}$, which is
stronger than property $\A_{\aleph_0}(1)$. One can use this result to show that
 $M_\ga(F^\infty)$ has property $\A_1(1)$.

Another consequence of Theorem 2.1 is the following.

    \proclaim{Theorem 2.9} For any cardinal $\ga$,
    the map $\Phi:M_\ga(F^\infty)/{M_\ga(J)}\to M_\ga(B(\N_J))$
    defined by
    $$
    \Phi([f_{\alpha\beta}]+M_\ga(J))=[P_{\N_J}f_{\alpha\beta}|_{\N_J}]
    $$
    is an isometry.
    \endproclaim
    \demo{Proof}
    It is enough to show that
    $$
    \text{\rm dist}([f_{\alpha\beta}],
    M_\ga(J))=\|[P_{U\N_J}U^*f_{i j}(S_1,\dots,S_n)U|_{U\N_J}]\|,
    $$
    where $U:\F^2\to\F^2$ is the flipping operator.
    For each $[g_{\alpha\beta}]\in M_\ga(J)$, we have
    $$
    \|[f_{\alpha\beta}+g_{\alpha\beta}]\|=
    \|[U^*(f_{\alpha\beta}+g_{\alpha\beta})U]\| 
    \ge
    \|[P_{U\N_J}U^*(f_{\alpha\beta}+g_{\alpha\beta})U|_{U\N_J}]\|.
    $$
    Since $g_{\alpha\beta}\in J$, according to Proposition 1.3, we have
    $P_{\N_J}g_{\alpha\beta}|_{\N_J}=0$. Since $U$ is an unitary operator with $U=U^*$, it is easy to
    see that
    $P_{U\N_J}U^*g_{\alpha\beta}U|_{U\N_J}=0$.
    Combining this with the above inequality, we obtain
    $$
    \text{dist\,}([f_{\alpha\beta}], M_\ga(J))\ge
    \|[P_{U\N_J}U^*f_{\alpha\beta}U|_{U\N_J}]\|.\tag 2.11
    $$
    It remains to prove the converse inequality.
    Since $U^* f_{\alpha\beta}U$ commutes with $S_1,\dots,S_n$, and $  U\N_J$ is invariant
    to $S_1^*,\dots,S_n^*, (U^*f_{\alpha\beta}U)^*$,  it is clear that
    $
    [P_{U\N_J}U^*f_{\alpha\beta}U|_{U\N_J}]
    $
    commutes with
   $P_{U\N_J}S_i|_{U\N_J}\otimes I_\ga$ for each $i=1,2,\dots,n$.
    We can apply Theorem 2.1 to find
    $[\psi_{\alpha\beta}]\in M_\ga(F^\infty)$ such that
    $\|[\psi_{\alpha\beta}]\|=\|[P_{U\N_J}U^*f_{\alpha\beta}U|_{U\N_J}] \|$
    and
    $[P_{U\N_J}U^*\psi_{\alpha\beta}U|_{U\N_J}]=
     [P_{U\N_J}U^*f_{\alpha\beta}U|_{U\N_J}]
     $.
     According to Proposition 1.3, we infer that
     $[\phi_{\alpha\beta}]:= [\psi_{\alpha\beta}-f_{\alpha\beta}]
     \in M_\ga(J)$.
     Therefore,
     $$
      \|[P_{U\N_J}U^*f_{\alpha\beta}U|_{U\N_J}] \|=
     \|[f_{\alpha\beta}+\phi_{\alpha\beta}]\|\ge \text{\rm dist}
     ([f_{\alpha\beta}], M_\ga(J)).
     $$
     Combining this with (2.11), we complete the proof.
$\blacksquare$
\enddemo

Recall that for each $k\geq1$, $M_k(F^\infty/J)=M_k(F^\infty)/M_k(J)$
(see [R]).  Hence, as an immediate consequence of Theorem 2.9, we obtain
the following.

\proclaim{Corollary 2.10} The map 
$\Phi: F^\infty/J\to P_{\N_J}F^\infty|_{\N_J}$ defined by $\Phi(f)= P_{\N_J}f|_{\N_J}$ is a completely isometric representation.
\endproclaim

Let $\P_m$ be the set of all polynomials in $\F^2$ of degree $\le m$, and denote
 \linebreak $\F^\infty_{\ge m}= \F^\infty \cap (\F^2 \ominus \P_m).$
Setting $J=\F^\infty_{\ge m}$ in Corollary 2.10, we can deduce the Carath\' eodory
interpolation theorem on Fock spaces [Po7].

\proclaim {Corollary 2.11} Let $p\in \P_m$ be fixed.
Then 
$${\text{\rm dist}}(p, \F^\infty_{\ge m})=\| P_{\P_m}p(S_1,\dots,S_n)|_{\P_m}\|.$$
\endproclaim

Let us remark that Theorem 2.9 is no longer true 
if we replace $F^\infty$ by the noncommutative disk algebra 
$\Al_n$ and $J$ is a closed, 2-sided ideal of $\Al_n$.
To see this, let $\lambda\in \Bbb C^n$ be of norm one and let
$J:=\{ \psi \in \Al_n : \psi(\lambda_1,\dots,\lambda_n)=0  \text{ and }\langle\psi,1\rangle=0 \}.$
It is easy to see that $\N_J$ is the span of $1$ and so is $\N_{J_w}$ (see also Example 3.6).  Then
$J_w= \{ \psi \in F^\infty : \langle\psi, 1\rangle=0 \}$.  
If one takes a polynomial $p\in \P$ such that $\langle p,1\rangle=0$ but
$p(\lambda_1,...,\lambda_n)\not = 0$, then
	$\text{\rm dist}(p,J)>0$ but $\text{\rm dist}(p,J_w)=0.$
Therefore,
$$
\text{\rm dist} (p,J)\neq \text{\rm dist}(p,J_w)=\| P_{\N_{J_w}}f|_{\N_{J_w}}\|=
\| P_{\N_J}f|_{\N_J}\|.
$$

However, we will show that $\Al_n/J$ is completely isometrically 
isomorphic to
$P_{\N_J}\Al_n|_{\N_J}$, for certain closed ideals $J$ of $\Al_n$.

\proclaim{Proposition 2.12} 
Let $\lambda_1,\dots,\lambda_k\in \Bbb B_n$ and define
$$\eqalign{
J &=\{\varphi\in \Al_n : \varphi(\lambda_j)=0 \hbox{ for every }
j=1,2,\dots,k\}, \hbox{ and}\cr
J_w &=\{\varphi\in F^\infty : \varphi(\lambda_j)=0 \hbox{ for every }
j=1,2,\dots,k\}.\cr}$$
Then the map
$\Psi:\Al_n/J\to P_{\N_J}\Al_n|_{\N_J}$ defined by
$\Psi(f+J)= P_{\N_J}f|_{\N_J}$
is a completely isometric representation.
\endproclaim
\demo{Proof}
According to Corollary 2.10 and Lemma 1.4, for any $f\in \Al_n$,
$$
\text{\rm dist}(f,J_w)=\|\ P_{\N_{J_w}}f|_{\N_{J_w}}\|=
\|\P_{\N_J}f|_{\N_J}\|.
$$	
Therefore, it is enough to prove that
$
\text{\rm dist}(f,J_w)=\text{\rm dist}(f,J).
$
Let us define  $\Phi:\Al_n/J\to F^\infty/J_w$ by 
$\Phi(\varphi+J)=\varphi+J_w$.  Notice that $\Phi$ is contractive.
We shall prove that for every $\varphi\in F^\infty$ with
$\|\varphi+J_w\|=1$, there exists
$\psi\in \Al_n$ such that $\|\psi+J\|=1$ and
$\Phi(\psi+J)=\varphi+J_w$.

  Assume that $\|\varphi\|=\|\varphi+J_w\|=1$ and find
$\varphi_k\in \Al_n$ such that $\|\varphi_k\|\leq 1$ and
$\varphi_k\to \varphi$ in the WOT.  Since $\Al_n/J$ is finite
dimensional, we assume (after passing to a subsequence) that
$\varphi_k+J$ converges to $\psi+J$ in the norm of $\Al_n/J$ for
some $\psi\in \Al_n$.
Then $\|\psi+J\|\leq 1$ and there exists a sequence $\eta_k\in J$
such that $\varphi_k+\eta_k\to \psi$ in the norm topology of $\Al_n$.
Then $\eta_k = (\eta_k+\varphi_k)-\varphi_k\to \psi-\varphi$ in the WOT
of $F^\infty$.
Since $\eta_k\in J\subset J_w$ for each $k$ and since $J_w$ is
WOT-closed, we have that 
$\psi-\varphi\in J_w$.  Therefore,
$\varphi+J_w=\psi+J_w=\Phi(\psi+J)$.  
 Since $\dim \Al_n/J=\dim F^\infty/J_w$, it is clear now that  $\Phi$ is isometric.
The argument works also when passing to matrices, so the map $\Phi$ is completely isometric.
$\blacksquare$
\enddemo

Combining Theorem 2.4 with Proposition 2.12, we infer the following Nevanlinna-Pick interpolation theorem for the noncommutative disc algebra $\Al_n$. For simplicity,
we consider only the scalar case.

\proclaim{Corollary 2.13}
Let  $\lambda_1,\dots,\lambda_k\in \Bbb B_n$, and  $w_1,\dots,w_k\in \C$. Then the matrix 
$$\left[{1-w_j{\bar w}_i \over 1-\langle\lambda_j,\lambda_i\rangle}
\right]_{i,j=1,\dots,k}$$
is positive definite if and only if for any $\epsilon>0$ there exists $f\in \Al_n$,
$\|f\|_\infty\leq1+\epsilon$, such that $f(\lambda_j)=w_j$ for every $j=1,\dots,k$.
\endproclaim

\bigskip

\leftline {\bf 3. Poisson Transforms and von Neumann Inequalities}

In [Po9], the second author found an elementary proof of the inequality (1) 
based on  noncommutative Poisson transforms associated to row contractions.   In this section, we will recall this construction (see [Po9, Section 8]) in a particular case and use it to obtain new results.

As in [Po1], $T=[T_1,\dots,T_n]$ is  called $C_0$-contraction if $T$ is a contraction and
$$\text{\rm SOT-}\lim_{k\to\infty}\sum_{\alpha\in \Bbb F_n^+,|\alpha|=k}T_\alpha T_\alpha^*=0.\leqno{(3.1)}$$
Recall that the sequence 
$\{\sum_{|\alpha|=k}T_\alpha T_\alpha^* : k\geq0\}$
of positive operators is non-increasing, and that   (3.1) holds
if and only if  \ $\sum_{|\alpha|=k}\|T_\alpha^*h\|^2\to 0$ for every $h\in \H$.

Suppose that $T=[T_1,\dots,T_n]$ is  a $C_0$-contraction and let  $\Delta:=(I_\H-\sum_{i=1}^n T_iT_i^*)^{1\over 2}$. Since $$\sum_{|\alpha|=k}T_\alpha\Delta^2 T_\alpha^*
=\sum_{|\alpha|=k}T_\alpha T_\alpha^*-\sum_{|\alpha|=k+1}T_\alpha T_\alpha^*,$$
it is clear that
 $\sum_{\alpha\in \Bbb F_n^+} T_\alpha \Delta^2 T_\alpha^*
=I_\H-\lim_{k\to\infty}\sum_{|\alpha|=k+1}T_\alpha T_\alpha^*=I_\H$.

The Poisson Kernel $K=K(T)$
associated to $T=[T_1, \dots, T_n]$ is the linear map
$$K:H\to \F^2\otimes \H\quad\quad\text{defined by}\quad\quad
Kh=\sum_{\alpha\in \Bbb F_n^+} e_\alpha\otimes \Delta T_\alpha^*h.$$
Since $\sum_\alpha T_\alpha \Delta^2 T_\alpha^* = I_\H$, 
$K$ is an isometry. 
 It is easy to check that, for each
$\alpha\in\Fr_n^+$, \   $(S_\alpha^*\otimes I)Kh=KT_\alpha^*h$. Hence,
 for every $\alpha,\beta\in \Fr_n^+$,
$$K^*\bigl[ S_\alpha S_\beta^*\otimes I\bigr]K=T_\alpha T_\beta^*.\tag 3.2$$
The map $\Psi:B(\F^2)\to B(\H)$ defined by $\Psi(A)=K^*[A\otimes I]K$
is clearly unital, completely contractive (hence, completely positive),
and $w^*$-continuous.
  Moreover, for each $\alpha,\beta\in\Fr_n^+$,
$\Psi(S_\alpha S_\beta^*)=T_\alpha T_\beta^*$.  
The restriction of $\Psi$ to $F^\infty$, which is denoted by $\Psi_T$,
provides a WOT-continuous
$F^\infty$-functional calculus for the $C_0$-contractions $T=[T_1,\dots,T_n]$, which is a particular case of [Po6].  That is,
$$\Psi_T:F^\infty \to B(\H)\quad\text{satisfies}\quad
\Psi_T(\varphi(S_1,\dots,S_n))=\varphi(T_1,\dots,T_n)\leqno{(3.3)}$$
for every $\varphi\in F^\infty$.

\smallbreak

Suppose now that $T=[T_1,\dots,T_n]$ is a row contraction.
For each $0<r<1$, let $K_r=K_r(T)$ be the Poisson Kernel associated to $[rT_1,\dots rT_n]$,
which is clearly a $C_0$-contraction.  Let $C^*(S_1,\dots,S_n)$ be the
$C^*$-algebra generated by $S_1,\dots,S_n$, the extension through compacts of the Cuntz algebra ${\Cal O}_n$ (see [Cu]). The Poisson Transform associated to $T=[T_1,\dots,T_n]$ is
the linear  map
$$\Phi_T:C^*(S_1,\dots,S_n)\to B(\H)\quad\text{defined by}\quad
\Phi_T(f)=\lim_{r\to 1} K_r^*[f\otimes I]K_r\leqno{(3.4)}$$
(in the uniform topology of $B(\H)$).
It is easy to see that $\Phi_T$ is
unital, completely contractive, and for every $\alpha,\beta\in\Fr_n^+$, 
$\Phi_T(S_\alpha S_\beta^*)=T_\alpha T_\beta^*$.  
Inequality (1) from the introduction follows by
restricting $\Phi_T$ to $\Al_n$.

\medbreak


A simple consequence of the noncommutative Poisson transform is the following result which turns out to be crucial for the rest of this paper.

 \proclaim{Proposition 3.1} Let $T=[T_1,\dots,T_n]$ be a  $C_0$-contraction
with its Poisson Kernel $K$, and let $\N$ be a subspace of $\F^2$ invariant 
under $S_1^*,\dots,S_n^*$.  If $K$ takes values in 
$\N\otimes\H$, then there exists a unital, completely contractive,
$w^*$-continuous map $\Phi:B(\N)\to B(\H)$ such that for every
$\alpha,\beta\in\Fr_n^+$, $\Phi(B_\alpha B_\beta^*)=T_\alpha T_\beta^*$,
where $B_k=P_\N S_k{}{\vert_\N}$ for any $ k=1,\dots, n$.
\endproclaim

\demo{Proof}  Since  $\N\subset\F^2$ is an invariant subspace
of  $S_1^*,\dots,S_n^*$, for every 
$\alpha,\beta\in\Fr_n^+$,  $P_\N S_\alpha S_\beta^* {}{\vert_\N}=B_\alpha B_\beta^*$.
By hypothesis, $K= (P_\N\otimes I)K$.  Hence, and according to (3.2), for each $\alpha,\beta\in \Bbb F_n^+$, we have
$$\eqalign{ T_\alpha T_\beta^*
&= K^*[S_\alpha S_\beta^*\otimes I]K=
	 K^*(P_\N\otimes I)[S_\alpha S_\beta^*\otimes I](P_\N\otimes I)K\cr
&=K^*[P_\N S_\alpha S_\beta^* P_\N\otimes I]K= 
K^*[B_\alpha B_\beta^*\otimes I]K.\cr}\leqno{(3.5)}$$
To complete the proof, define $\Phi:B(\N)\to B(\H)$ 
by $\Phi(A)=K^*[A\otimes I]K$. $\blacksquare$
\enddemo
\medbreak

\proclaim{Remark 3.2}  If  
$T=[T_1,\dots,T_n]$ is a  contraction and  its Poisson kernel $K_r$ takes values in $\N\otimes\H$ for every
$0<r<1$, 
then there is  a unital, 
completely contractive map $\Phi:C^*(B_1,\dots,B_n)\to B(\H)$
satisfying $\Phi(B_\alpha B_\beta^*)=T_\alpha T_\beta^*$ for all
$\alpha,\beta\in\Fr_n^+$.
\endproclaim
\demo{Proof}
  It follows from (3.5)
that 
$$\lim_{r\to 1} K_r^*[B_\alpha B_\beta\otimes I]K_r=\lim_{r\to 1}
r^{|\alpha|}T_\alpha r^{|\beta|}T_\beta^*=T_\alpha T_\beta^*.$$  Hence,
the map $B_\alpha B_\beta^*\mapsto T_\alpha T_\beta^*$,
defined on span$\{B_\alpha B_\beta^*:\alpha,\beta\in\Fr_n^+\}$, is 
completely contractive.  By [Ar1], it can be extended to a unital, 
completely contractive map $\Phi:C^*(B_1,\dots,B_n)\to B(\H)$
satisfying $\Phi(B_\alpha B_\beta^*)=T_\alpha T_\beta^*$ for all
$\alpha,\beta\in\Fr_n^+$. $\blacksquare$

\enddemo

 \smallbreak

To illustrate Proposition 3.1 and Remark 3.2, we will consider a row contraction
 $T=[T_1,\dots,T_n]$ satisfying the following commutation relations
$$T_jT_i=\lambda_{ji}T_iT_j\quad\text{for every}\quad 1\leq i<j\leq n,
\leqno{(3.6)}$$
where $\lambda_{ij}\in\C$ for $1\leq i<j\leq n$.

\proclaim{Example 3.3} There
exists a subspace $\N=\N(\{\lambda_{ij}\})$ of $\F^2$,
invariant under $S_1^*,\dots,S_n^*$, such that the
operators $B_k=P_\N S_k{}{\vert_\N}$, $k=1,\dots, n$, satisfy $(3.6)$ and 
for every  row contraction $T=[T_1,\dots,T_n]$
satisfying $(3.6)$, there exists 
a unital completely contractive linear  map $\Phi:C^*(B_1,\dots,B_n)\to B(\H)$
such that  
$\Phi(B_\alpha B_\beta^*)=T_\alpha T_\beta^*$ for any $\alpha,\beta\in\Fr_n^+$.
\endproclaim

\demo{Proof}  Fix $k\in\Na$, and consider
$\alpha=g_{i_1}g_{i_2}\cdots g_{i_k}\in \Fr_n^+$ satisfying
$i_1\leq i_2\leq\cdots\leq i_{k}$, and a permutation  
$\pi\in\Pi_k$ on $\{1,2,\dots,k\}$.  Then, from (3.6), 
$$T_{\pi(\alpha)}=\epsilon_{\pi(\alpha)}T_\alpha,\quad\text{where}\quad
\epsilon_{\pi(\alpha)}:=\prod_{{j<\ell \atop {\pi(j)}> {\pi(\ell)}}}
\lambda_{i_{\pi(j)} i_{\pi(\ell)}}$$
and $\pi(\alpha):=g_{i_{\pi(1)}}g_{i_{\pi(2)}}\cdots g_{i_{\pi(k)}}$.
 Let $\N(\{\lambda_{ij}\})$ be the subspace of $\F^2$ defined by
$$\N(\{\lambda_{ij}\}):=\overline{\text{span}}\biggl\{
\sum_{\pi\in\Pi_k} \overline{\epsilon}_{\pi(\alpha)} e_{\pi(\alpha)}:
\alpha=g_{i_1}g_{i_2}\cdots g_{i_k}\in \Fr_n^+,
i_1\leq\cdots\leq i_k, k\in\Na \biggr\}.$$
It is easy to see that if $T=[T_1,\cdots,T_n]$ is  a $C_0$-contraction
and satisfies (3.6), then its Poisson kernel takes values in
$\N(\{\lambda_{ij}\})\otimes\H$.
  Indeed, 
$$Kh=\sum_{k=0}^\infty
\sum_{\alpha=g_{i_1}\cdots g_{i_k}\atop i_1\leq\cdots\leq i_k}
\sum_{\pi\in\Pi_k} e_{\pi(\alpha)}\otimes \Delta T_{\pi(\alpha)}^*h=
\sum_{k=0}^\infty
\sum_{\alpha=g_{i_1}\cdots g_{i_k}\atop i_1\leq\cdots\leq i_k}
v_\alpha \otimes\Delta T_\alpha^*h,$$
where $v_\alpha=\sum_{\pi\in\Pi_k} \bar{\epsilon}_{\pi(\alpha)}e_{\pi(\alpha)}
\in \N(\{\lambda_{ij}\}).$  One can verify directly, 
from the definition of $\N(\{\lambda_{ij}\})$, that this space is
invariant under $S_1^*,\dots,S_n^*$, although it
is easier to check that $\N(\{\lambda_{ij}\})=\N_J$, where
$J$ is the WOT-closed, 2-sided ideal in $\F^\infty$ generated by 
$\{e_j\otimes e_i-\lambda_{ji}e_i\otimes e_j:
1\leq i<j\leq n\}$. 
Then, from Proposition 1.3,  the $B_k$'s satisfy (3.6).  
The rest of the statement of Example 3.3 is an immediate consequence of 
 Remark 3.2. $\blacksquare$
\enddemo

The case where $\lambda_{ji}=1$ for $1\leq i< j\leq n$
appears in [Ath] and [Ar2].  In this situation, condition (3.6) means that the
$T_i$'s are commuting and $\N(\{\lambda_{ji}\})$ is the symmetric
Fock space.  If $\lambda_{ji}=-1$ for $1\leq i<j\leq n$, then the  $T_i$'s are anti-commuting and
$\N(\{\lambda_{ji}\})$ is the anti-symmetric Fock space.

\proclaim{Example 3.4}If $J_k$ is a  \text{\rm WOT}-closed, $2$-sided ideal generated by 
some elements in
$\text{\rm span}\{e_\alpha: \ |\alpha|=k\}$, then a similar result to Example $3.3$ holds for any contraction $T=[T_1,\dots,T_n]$ such that $\phi(T_1,\dots,T_n)=0$ for each 
$\phi\in J_k$.
\endproclaim
In Section 4, we will consider the $F^\infty$--functional calculus associated to row contractions satisfying (3.6), or as in Example 3.4.

Let $\phi\in F^\infty$ and let $J_\phi$ be the WOT-closed, $2$-sided ideal generated by $\phi$ in $F^\infty$. If $\N_{J_\phi}\neq \{0\}$, then there is a nontrivial $C_0$-contraction  $T=[T_1,\dots,T_n]$ such that $\phi(T_1,\dots,T_n)=0$. 
Indeed, define $T_i:=P_{\N_{J_\phi}} S_i|_{\N_{J_\phi}}$. 
According to [Po1], it is clear that 
$T=[T_1,\dots,T_n]$ is a $C_0$-contraction.
Since the $F^\infty$- functional calculus associated to $C_0$-contractions is WOT- continuous.  It is easy to see that 
$\phi(T_1,\dots,T_n)=P_{\N_{J_\phi}} \phi(S_1,\dots,S_n)|_{\N_{J_\phi}}=0$
(see also Lemma 4.4).

\proclaim{Lemma 3.5}  Suppose that 
$T=[T_1,\dots,T_n]$ is a $C_0$-contraction 
with  its Poisson Kernel $K$, and that $J$ is a \text{\rm WOT}-closed, $2$-sided ideal of 
$F^\infty$ such that for every $\varphi\in J$, $\varphi(T_1,\dots,T_n)=0$.
Then $K$ takes values in $\N_J\otimes \H$. Consequently, $\N_J\not=(0)$.
\endproclaim

\demo{Proof}  
For any polynomial $p\in \P,\ p=\sum_\alpha a_\alpha e_\alpha$, we have
$$
\align
\langle Kk,p\otimes h \rangle 
&=\sum_\alpha \bar{a}_\alpha \langle k, T_\alpha \Delta h\rangle
=\biggl\langle k,\biggl(\sum_\alpha a_\alpha T_\alpha\biggr)\Delta h
\biggr\rangle\\
&=\langle k,p(T_1,\dots,T_n)\Delta h\rangle
\endalign
$$
for any $h,k\in\H$.
 Since the $F^\infty$-functional calculus for $C_0$-contractions is  WOT-continuous and $\P$ is WOT-dense in $F^\infty$ we deduce that for any $\varphi\in J$ and $h, k\in\H$,
$$
\langle Kk,\varphi\otimes h \rangle 
=
\langle k, \varphi(T_1,\dots,T_n)\Delta h\rangle=0.
$$
Since $\M_J$ is the closure of $J$ in $\F^2$, we see that
for every $k\in\H$, 
$$Kk\in \bigl(\M_J\otimes H\bigr)^\perp=\N_J\otimes H.$$
This completes the proof. 
$\blacksquare$
\enddemo

If $T=[T_1,\dots,T_n]$ is a $C_0$-contraction, then
$$J_1:=\{\varphi\in F^\infty: \varphi(T_1,\dots,T_n)=0\}=\text{Ker\,}\Psi_T$$ is
a WOT-closed, 2-sided ideal of $F^\infty$. Similarly,
$$J_2:=\{\varphi\in \Al_n: \varphi(T_1,\dots,T_n)=0\}=\text{Ker\,}\Phi_T$$ is
a closed, 2-sided ideal of $\Al_n$.  Lemma 3.5 is stated for $F^\infty$,
but it holds true also for $\Al_n$. Therefore  $\N_{J_1}\not=\{0\}$  and $\N_{J_2}\not=\{0\}$.
 Let us remark that if $[T_1,\dots,T_n]$ is just a  contraction 
(not necessarily $C_0$), then $\N_{J_2}$ may be zero.

\medbreak

\noindent{\bf Example 3.6.} (Point evaluations) 
Let $\lambda_i\in\C$, $i=1,\dots, n$, be such that \linebreak $\sum_{i=1}^n|\lambda_i|^2<1$.
Then  $\lambda=[\lambda_1,\dots,\lambda_n]$ is a $C_0$-contraction,
and hence, $J_w:=\{\varphi\in F^\infty:\varphi(\lambda_1,\dots,
\lambda_n)=0\}$ is a WOT-closed, 2-sided ideal of $F^\infty$.  It is known
that $\N_{J_\lambda}=\text{span}\{z_\lambda\}$ where
$z_\lambda=1+\sum_{k\geq1}(\lambda_1e_1+\cdots+\lambda_ne_n)^{\otimes k}$
and
$\varphi(\lambda_1,\dots,\lambda_n)=\langle\varphi,z_\lambda\rangle$
for every $\varphi\in F^\infty$
(see [APo], [Ar2], and [DP1]).  Notice that if $\sum_{i=1}^ n|\lambda_i|^2=1$, then
$J=\{\varphi\in\Al_n:\varphi(\lambda_1,\dots,\lambda_n)=0\}$
is a closed, 2-sided ideal of $\Al_n$ but one can check that
$\N_J=\{0\}$.

\medbreak

Combining Proposition 3.1 and Lemma 3.5, we obtain the following.

\proclaim{Theorem 3.7} Let $T=[T_1,\dots,T_n]$ be a $C_0$-contraction, 
and let $J$ be a \text{\rm WOT}-closed, $2$-sided ideal of $F^\infty$
such that for every $\varphi\in J$, $\varphi(T_1,\dots,T_n)=0$,
then there exists a unital, completely contractive,
$w^*$-continuous map $\Phi:B(\N_J)\to B(\H)$ such that for every
$\alpha,\beta\in\Fr_n^+$, $\Phi(B_\alpha B_\beta^*)=T_\alpha T_\beta^*$,
where $B_k=P_{\N_J} S_k{}\vert_{N_J}, k=1,\dots, n$.
\endproclaim

One can easily see
that there is an $\Al_n$-version of this theorem corresponding to closed, $2$-sided ideals in $\Al_n$, $J\subset \text{Ker\,}\Phi_T$ with $\N_J\neq \{0\}$.
Let  $T=[T_1,\dots,T_n]$ 
be  a contraction, and let $J\subset \text{Ker\,}\Phi_T$ be a closed, $2$-sided ideal of $\Al_n$ such that $\N_J\neq \{0\}$. Notice that Remark 3.2 holds true if we take $\N=\N_J$.

Given $T=[T_1,\dots,T_n]$ a  $C_0$-contraction with Poisson kernel $K$,
the best von Neumann inequality given by Proposition 3.1 comes from the smallest
subspace $\N_T$ of $\F^2$ which is invariant under $S_1^*,\dots,S_n^*$
and such that $K$ takes values in $\N_T\otimes \H$.  
It is not hard to see that
$\N_T=\overline{\text{span}}\, \{
\langle T_\alpha^* k,\Delta h\rangle e_\alpha :h,k\in\H; \alpha\in\Fr_n^+ \}.$
First notice that $\N_T$ is the smallest $\N$ such that $K$ takes values
in $\N\otimes \F^2$, and then notice that $\N_T$ is invariant under
$S_1^*,\dots,S_n^*$.

\bigskip
\leftline {\bf 4. $\W(B_1,\cdots,B_n)$ and 
$F^\infty/J$-functional calculus for row contractions}

In this section $J$ will be a $w^*$-closed, 2-sided ideal of $F^\infty$.
Recall that $\N_J$ is the orthogonal complement of the image of $J$
in $\F^2$ and that $B_k=P_{\N_J}S_k{}\vert_{\N_J}$ for $k=1,\dots, n$.
We define $\W(B_1,\dots,B_n)$ to be $w^*$-closure of of the algebra
generated by the $B_k$'s and the identity.

We will prove that $F^\infty/J$ is canonically isomorphic to
$\W(B_1,\dots,B_n)$.  
We will describe the commutant of $\W(B_1,\dots,B_n)$ and
will show that $\W(B_1,\dots,B_n)$ is the double commutant of
$\{B_1,\dots,B_n\}$.
We will show that $\W(B_1,\dots,B_n)$ has
the $\A_1(1)$ property and hence the $w^*$ and WOT topologies agree
on this algebra.
Finally, we will develop a 
$F^\infty/J$-functional calculus for row contractions.  

A direct consequence of 
Proposition 1.2 and Corollary 2.10 is the following.

\proclaim{Theorem 4.1}  The map
$\Psi:F^\infty/J\to B(\N_J)$ defined by
$$\Psi(\varphi+J)=P_{\N_J}\varphi(S_1,\dots,S_n){}\vert_{\N_J}$$
is a completely isometric isomorphism onto $P_{\N_J} F^\infty|_{\N_J}$, and a homeomorphism relative to the $w^*$- topology on $F^\infty/J$ and the \text{\rm WOT}-topology on $P_{\N_J} F^\infty|_{\N_J}$.
  \endproclaim

\demo{Proof}  
Since the fact that $\Psi$ is a completely isometric homomorphism was already 
proved in Corollary 2.10 (see also Section 5), we only have to 
prove that $\Psi$ is a $w^*$-WOT homeomorphism.

By Proposition 1.2, $\varphi_i+J \to \varphi+J$ in the $w^*$ topology
iff for every $\xi_1,\xi_2\in\N_J$, $\langle \varphi_i \otimes\xi_1,\xi_2\rangle\to
\langle \varphi \otimes\xi_1,\xi_2\rangle$.  This is clearly equivalent to
$P_{\N_J}\varphi_i{}\vert_{\N_J}\to P_{\N_J}\varphi{}\vert_{\N_J}$
in the weak operator topology.
$\blacksquare$
\enddemo
\medbreak

Using  again the noncommutative commutant lifting theorem [Po3], we can prove the following.

\proclaim{Proposition 4.2}The algebra $P_{\N_J}F^\infty|_{\N_J}$ is the \text{\rm WOT}-closed algebra generated by $P_{\N_J} {S_i}\vert_{\N_J},\ i=1,\dots,n,$
and the identity. Moreover, we have
$$
P_{\N_J}F^\infty|_{\N_J}=\{ P_{\N_J}U^*F^\infty U|_{\N_J}\}'=  \{P_{\N_J}F^\infty|_{\N_J}\}''.\tag 4.1
$$
\endproclaim
\demo{Proof}
We first show that $P_{\N_J}F^\infty|_{\N_J}$ is weakly closed.
Notice that $\N_J$ is an invariant subspace of $U^* S_i^* U$ for each $i=1,\dots,n$, and $[P_{\N_J} {S_1}\vert_{\N_J},\dots, P_{\N_J} {S_n}\vert _{\N_J}]$ is a $C_0$-contraction with the minimal isometric dilation $[U^*S_1U,\dots, U^*S_n U]$. According to the commutant lifting theorem,   
$X\in\{ P_{\N_J}U^*S_iU|_{\N_J}:\ i=1,\dots,n\}'$ if and only if $X= P_{\N_J}Y|_{\N_J}$ for some $Y\in \{U^*S_1U,\dots, U^*S_n U\}'$. Using
[Po7], we get $Y=f(S_1,\dots, S_n),\ f\in \F^\infty$.
Now, it is clear that $P_{\N_J}F^\infty|_{\N_J}=\{ P_{\N_J}U^*F^\infty U|_{\N_J}\}'$ and, hence,  $P_{\N_J}F^\infty|_{\N_J}$ is a WOT-closed algebra.

 Since the polynomials in $S_1,\dots, S_n$ are WOT-dense in 
$F^\infty$, it is clear that $P_{\N_J}F^\infty|_{\N_J}$ is the WOT-closed algebra generated by $P_{\N_J} {S_i}\vert_{\N_J},\ i=1,\dots,n,$
and the identity. The second equality in (4.1) follows in a similar manner.
 $\blacksquare$
\enddemo

\proclaim{Proposition 4.3} The algebra  $\W(B_1,\dots,B_n)$ has 
property $\A_1(1)$. Moreover, \linebreak
$\W(B_1,\dots,B_n)=P_{\N_J}F^\infty|_{\N_J}$.
\endproclaim

\demo{Proof}  Let $f\in \W(B_1,\dots,B_n)_*$ satisfy $\|f\|=1$.
Since $$\W(B_1,\dots,B_n)_*\equiv c_1(\N_J)/{}^\perp\W(B_1,\dots,B_n),$$
for each $\epsilon>0$,
find $g\in c_1(\N_J)$ satisfying $\|g\|\leq 1+\epsilon$ and
$f=g+{}^\perp\W(B_1,\dots,B_n)$.  
Let $i_{\N_J}:\N_J\to \F^2$ be the inclussion and notice that
$(i_{\N_J})^*=P_{\N_J}$.  It is easy to check that 
$i_{\N_J}\circ g\circ P_{\N_J}\in {}^\perp J\subset c_1(\F^2)$.
Then, by Proposition 1.2, there exists $\varphi,\psi\in\N_J$ satisfying
$\|\varphi\|_2\|\psi\|_2\leq(1+\epsilon)\|g\|
\leq(1+\epsilon)^2$ such that, for every $\eta\in F^\infty$,
$\langle i_{\N_J}\circ g\circ P_{\N_J},\eta\rangle = \langle\eta\otimes\varphi_1,\varphi_2\rangle$.
Now, for each noncommutative polynomial $p\in \P$, we have
$$\eqalign{\langle f, p(B_1,\dots,B_n)\rangle
&=\langle g,p(B_1,\dots,B_n)\rangle 
		=\langle g, P_{\N_J}p(S_1,\dots,S_n){}\vert_{\N_J}\rangle\cr
&=\langle i_{\N_J}\circ g\circ P_{\N_J}, p(S_1,\dots,S_n)\rangle 
		=\langle p\otimes\varphi,\psi \rangle\cr
&=\langle p(S_1,\dots,S_n) \varphi, P_{\N_J}\psi \rangle
		=\langle P_{\N_J} p(S_1,\dots,S_n){} \vert_{\N_J}\varphi,\psi \rangle\cr
&=\langle p(B_1,\dots,B_n)\varphi,\psi \rangle.\cr}$$
Since $f$ is $w^*$-continuous, we prove the $\A_1(1)$ property.
The last part of the theorem follows from Proposition 4.2. $\blacksquare$
\enddemo

\medbreak

Let $J$ be the $w^*$-closed, 2-sided ideal of $F^\infty$ generated by
$S_2, S_3,\dots, S_n$.  It is easy to see that $\N_J$ is the closed
span of $e_1^{\otimes k}$ for $k\geq 0$,  and that $B_2=\cdots=B_n=0$.
Hence, $\W(B_1,\dots,B_n)=\W(B_1)$, where 
$B_1 e_1^{\otimes k}=e_1^{k+1}$.  Since $B_1$ is a unilateral
shift of multiplicity one, we use Proposition 4.3 to give an alternative
proof of the well known fact that $\W(B_1)$ has property $\A_1(1)$.
Moreover, it is also known (see [BFP, Theorem 4.16])
that $\W(B_1)$ does not even satisfy property $\A_2$.  In that sense, 
Proposition 4.3 is best possible.

Let us recall from [Po1] that a contraction $~[T_1,\dots,T_n]~$ is called
completely non-coisometric (c.n.c.) if there is no $~h\in\H,\ h\not=0~$
such that
$$
\sum\limits_{|\alpha|=k}\|T_f^*h\|^2=\|h\|^2,\quad\text{for any}
\ k\in\{1,2,\dots\}.
$$

Let $T=[T_1,\dots,T_n]$ be a c.n.c. contraction
and let $$\Psi_T:F^\infty\to B(\H),\quad \Psi_T(f)=f(T_1,\dots,T_n),$$
be the $F^\infty$-functional calculus associated to $T$.
In this section, we prove that if $J$ is a WOT-closed, $2$-sided ideal of $F^\infty$ with $J\subset \text{\rm Ker\,}\Psi_T$, then there is a WOT-continuous,
$F^\infty/J$--functional calculus associated to $T$.

\proclaim{Lemma 4.4} Let $B=[B_1,\dots,B_n]$ and let $\Psi_B$ be the $F^\infty$-functional calculus associated to it. Then 
$$\W(B_1,\dots,B_n)=\Psi_B(F^\infty)=\{f(B_1,\dots,B_n): \ f\in F^\infty\}.
$$
\endproclaim
\demo{Proof}
 According to Proposition 4.3, it is enough to prove that
   $$
   f(B_1,\dots,B_n)=P_{\N_J}f(S_1,\dots,S_n)|_{\N_J}\tag 4.2
   $$
   for any $f\in \F^\infty$. Since $B_i=P_{\N_J}S_i|_{\N_J}$, (4.2) holds for
   polynomials, and consequently for elements in the noncommutative disc algebra
   $\Al_n$. Since $B=[B_1,\dots,B_n]$ is a $C_0$-contraction, according to the
    $F^\infty$-functional calculus, we have
    $$
    \align
     f(B_1,\dots,B_n):&=\text{SOT-}\lim_{r\to 1}f_r(B_1,\dots,B_n) \\
     &=\text{SOT-}\lim_{r\to 1} P_{\N_J}f_r(S_1,\dots,S_n)|_{\N_J}
     =P_{\N_J}f(S_1,\dots,S_n)|_{\N_J}
     \endalign
     $$
for any $f\in F^\infty$.
$\blacksquare$
\enddemo

\proclaim{Theorem 4.5} Let $T=[T_1,\dots,T_n]$ be a c.n.c. contraction
and let $$\Psi_T:F^\infty\to B(\H),\quad \Psi_T(f)=f(T_1,\dots,T_n),$$
be the $F^\infty$-functional calculus associated to $T$. 
If $J$ is a \text{\rm WOT}-closed, $2$-sided ideal of $F^\infty$ 
with $J\subset \text{\rm Ker}\Psi_T$, then the map 
$$\Psi_{T, J}:\W(B_1,\dots B_n)\to B(\H);\quad \Psi_{T,J}(f(B_1,\dots,B_n)):=f(T_1,\dots,T_n),\tag 4.3
$$
is a \text{\rm WOT}-continuous, completely contractive homomorphism. 
In particular, for any $f\in F^\infty$,
$$
\|f(T_1,\dots,T_n)\|\le \|f(B_1,\dots,B_n)\|=\text{\rm dist}(f,J).
$$
\endproclaim
 \demo{Proof}
 We prove first that $\Psi_{T,J}$ is WOT-continuous.
  Let $f_i, f\in\F^\infty$ with
   $$\text{\rm WOT-}\lim_i f_i(B_1,\dots,B_n)=f(B_1,\dots,B_n).$$
  According to Lemma 4.4, we infer that
  $\text{\rm WOT-}\lim_i P_{\N_J}f_i|_{\N_J}= P_{\N_J}f|_{\N_J}$.
  Applying Proposition 4.1, we infer that
  $$ w^*{\text -}\lim_i(f_i +J)=f+J. \tag 4.4
  $$
  For each $h,k\in \H$,
  define $\Phi(f):=\langle \Psi_T(f)h,k\rangle$.
  Since $\Psi_T$ is WOT-continuous, $\Phi$ is WOT-continuous, and hence
  $w^*$-continuous.
  On the other hand, $\Psi(J)=0$, so that $\Phi\in{}^\perp J$.
  Since (4.4) holds, we deduce that $\lim_i \Phi(f_i)=\Phi(f)$, which is
   equivalent to
  $$\lim_i \langle f_i(T_1,\dots,T_n)h,k\rangle=\langle f(T_1,\dots,T_n)h,k\rangle
  $$
  for any $h,k\in \H$.

  According to the von Neumann inequality [Po5], for any $\psi\in J\subset \text{\rm Ker\,}\Psi_T$, we have
  $$
  \|f(T_1,\dots, T_n)\|=\|(f+\psi)(T_1,\dots, T_n)\|\leq \|f+\psi\|_\infty.
  $$
  Using Theorem 4.1, we infer that
  $$
  \align
  \|f(T_1,\dots, T_n)\|&\leq \text{\rm dist}(f,J)=\|P_{\N_J}f|_{\N_J}\|\\
  &=\|f(B_1,\dots, B_n)\|.
  \endalign
  $$
  In a similar manner, one can prove that $\Psi_{T,J}$ is a completely contractive
homomorphism. This completes the proof.
$\blacksquare$
\enddemo

 The following  $F^\infty$-extension is  related to Example 3.3.

\proclaim{Corollary 4.6}  Let $T=[T_1,\dots,T_n]$ be a c.n.c. contraction
 satisfying the following commutation relations
$$T_jT_i=\lambda_{ji}T_iT_j\quad\text{for every}\quad 1\leq i< j\leq n,
$$
where $\lambda_{ij}\in\C$ for $1\leq i<j\leq n$. If
$J$ is the \text{\rm WOT}-closed, $2$-sided ideal generated by 
$\{e_j\otimes e_i-\lambda_{ji}e_i\otimes e_j:
1\leq i< j\leq n\}$ in $\F^\infty$, then there is a \text{\rm WOT}-continuous functional calculus
given by $(4.3)$. 
\endproclaim

\bigskip

\leftline {\bf 5. Representations of Quotients of Dual Algebras}

Recall that an operator algebra is a closed subalgebra of $B(\H)$ and
that a dual algebra is a unital $w^*$-closed subalgebra of $B(\H)$.
In the late 60's, Cole (see [BD, pages 270--273]) proved that 
quotients of uniform algebras are operator algebras.  Shortly after,
Lummer and Bernard proved that quotients of operator algebras are isometrically
isomorphic to operator algebras.
In [Pi, Chapter 4] Pisier noted that these methods also show 
that quotients of operator algebras are completely isometrically isomorphic
to operator algebras.   In this Section we will follow these ideas
closely to obtain simple representations of quotients of dual algebras.
As an application, we give an alternative proof of Corollary 2.10 that 
does not depend on the commutant lifting theorem of [Po6].

\proclaim{Proposition 5.1}  Let $A$ be a unital, $w^*$-closed subalgebra of
the bounded operators on a separable Hilbert space $\H$
such that for each $k\geq 1$, $M_k(A)$ has property $\A_1(1)$,
and let $J$ be a $w^*$-closed, $2$-sided ideal of $A$.  Then there exists
a subspace $\E\subset \ell_2\otimes  \H$ such that the map
$\widehat{\Psi}:A/J\to B(\E)$ defined by $\widehat{\Psi}(a+J)
=P_\E (I_{\ell_2}\otimes a) {} {\vert_\E}$ is a 
completely isometric representation.
\endproclaim

\demo{Proof}  Let $x\in M_k(A)/M_k(J)$, $\|x\|=1$.  We claim that
for every $\epsilon>0$, there exists a subspace $E\subset\ell_2^k(\H)$
such that the map
$$\Psi_x:A/J\to B(E)\quad\text{defined by}\quad
\Psi_x(a+J)=P_E(I_{M_k}\otimes a){} {\vert}_E
\leqno{(5.1)}$$
is a completely contractive homomorphism which satisfies
$\|(I_{M_k}\otimes\Psi_x)(x)\|\geq 1-\epsilon$.  If we take direct sums 
$\oplus_{x,\epsilon>0}\Psi_x$, where $x$ runs over the unit ball of 
$M_k(A)/M_k(J)$, $k\geq1$, and $\epsilon>0$, we get a completely
isometric embedding of $A/J$.  It will be clear from the construction
that it is enough to take countably maps $\Psi_x$, so the proposition  follows.

Let $\epsilon>0$ and write $x=y+M_k(J)$, where $y=(y_{ij})\in M_k(A)$.
Find $f\in (M_k(A)/M_k(J))^*=M_k(J)^\perp$, $\|f\|=1$, such that
$\langle x,f\rangle =1$.  Since $({}^\perp M_k(J))^{**}=M_k(J)^\perp$,
we can find $g\in {}^\perp M_k(J)$, $\|g\|\leq1$, such that
$|\langle y,g\rangle -\langle y,f\rangle|<\epsilon$.  
Then $|\langle y,g\rangle |\geq 1-\epsilon$.  Since $M_k(A)$ has the
$\A_1(1)$ property, find $\varphi,\psi\in\ell_2^k(\H)$, \linebreak
$\|\varphi\|_2=\|\psi\|_2=1$ such that for each $\eta\in M_k(A)$,
$\langle g, \eta \rangle=\langle \eta \varphi, \psi\rangle$.

Let $E_1=\overline{\text{span}}\{ \eta \varphi : \eta \in M_k(A)\}\subset 
\ell_2^k(\H)$, $E_2=\overline{\text{span}}
\{ \xi \varphi : \xi \in M_k(J)\}
\subset E_1$, and $E=E_1\ominus E_2$.
Since $E_1$ and $E_2$ are invariant under $M_k(A)$,  
the map
$\Phi_E:M_k(A)\to B(\ell_2^k(\H))$,
defined by $\Phi_E( \eta )=P_E \eta  {} {\vert_E}$, 
is a completely contractive homomorphism
that vanishes on $M_k(J)$.  Hence, the map 
$\Psi_x( a +J)=\Phi_E(I_{M_k}\otimes a)$ of (5.1) is a well
defined completely contractive representation.  

Since $\varphi\in E_1$, $\psi\in E_2^\perp$, and $E_2^\perp$
is invariant under the adjoints of $M_k(A)$, it is easy to check that
$\langle g, y \rangle =\langle y\varphi,\psi\rangle=
\langle y P_E(\varphi),P_E(\psi)\rangle$.  Hence,
$$\bigl|\langle \Phi_E(y) P_E\varphi, P_E\psi\rangle\bigr|=
\bigl|\langle g,y\rangle\bigr|\geq 1-\epsilon.$$

For each $i,j\leq k$, let $E_{ij}=\Phi_E(e_{ij}\otimes 1_{A})$.
The $E_{ij}$'s are matrix units on $E$, which decompose
$E=\H_1\oplus \H_2\oplus \cdots \oplus \H_n$.  $E_{ii}$ is the
orthogonal projection onto $\H_i$ and $E_{ij}$ is a partial isometry 
from $\H_j$ onto $\H_i$.  
Note that $E_{ij}$ commutes with the range of $\Psi_x$,
$E_{ij}=E_{ii}E_{i1}E_{1j}E_{jj}$, and
$\Phi_E(e_{ij}\otimes y_{ij})=
\Phi_E((e_{ij}\otimes 1_A)(I_{M_k}\otimes y_{ij}))=
\Phi_E(e_{ij}\otimes 1_A)\Phi_E(I_{M_k}\otimes y_{ij})=E_{ij}\Psi_x(y_{ij}+J)$.  
Let $\varphi_j=E_{jj}P_E\varphi$ and $\psi_i=E_{ii}P_E\psi$. Then
$$\eqalign{
\langle\Phi_E(y)P_E\varphi,P_E\psi\rangle
&=\sum_{i,j\leq k} \langle\Phi_E(e_{ij}\otimes y_{ij})P_E\varphi,P_E\psi\rangle\cr
&=\sum_{i,j\leq k} \langle E_{ij}\Psi_x(x_{ij})P_E\varphi,P_E\psi\rangle\cr
&=\sum_{i,j\leq k} \langle\Psi_x(x_{ij})E_{1j}\varphi_j,E_{1i}\psi_i\rangle\cr
&=\langle (I_{M_k}\otimes \Psi_x)(x) \hat{\varphi},\hat{\psi}\rangle,\cr}$$
where $\hat{\varphi}=(E_{11}\varphi_1,\dots,E_{1k}\varphi_k)
\in\ell_2^k(E)$, 
$\hat{\psi}=(E_{11}\psi_1,\dots,E_{1k}\psi_k)
\in\ell_2^k(E)$.  Since
$\|\hat{\varphi}\|_2^2=\sum_{j\leq k}\|E_{ij}\varphi_j\|_2^2\leq
\sum_{j\leq k}\|\varphi_j\|_2^2=\|P_E\varphi\|_2^2\leq 1$, and,
similarly, $\|\hat{\psi}\|_2\leq 1$, we get
$$\bigl\|(I_k\otimes\Psi_x)(x)\bigr\|\geq
\bigl|\langle (I_k\otimes\Psi_x)(x)\hat{\varphi},\hat{\psi}\rangle\bigr|
=\bigl|\langle \Phi_E(y) P_E\varphi, P_E\psi\rangle\bigr|\geq 1-\epsilon,$$
which proves the claim.  
Finally, notice that the map $\Psi_x$ was determined by \linebreak $g\in {}^\perp M_k(J)
\subset (M_k(A))_*$.  Since $(M_k(A))_*$ is separable
for each $k\geq1$, it is enough to take only countably many maps.
$\blacksquare$
\enddemo

The proof of the next corollary follows easily from the proof of Proposition
5.1. Notice that property $\A_1(1)$ can be used to give 
more explicit representations
of quotient algebras than those appearing
in Theorem 3.2 of [CW] and Theorem 0.3 of [Mc].

\proclaim{Corollary 5.2} Let $A$ be a unital, $w^*$-closed subalgebra
of $B(\H)$ with the $\A_1(1)$ property and let $J\subset A$ be a $w^*$-closed
$2$-sided ideal.  Then for every $T\in A$,
$$\text{\rm dist}(T,J)=\sup\biggl\{\bigl\| P_{E_\varphi}T\vert_{E_\varphi}
\bigr\|_{B(\H)}:\varphi\in\H\biggr\},$$
where $E_\varphi=\overline{\text{\rm span}}\{a\varphi :a\in A\}\ominus
\overline{\text{\rm span}}\{b\varphi :b\in J\}\subset \H.$
\endproclaim

Moreover, it is well known that if $A\subset B(\H)$ is a unital
$w^*$-closed subalgebra of $B(\H)$, the ampliation 
$A^{(\infty)}=\{I_{\ell_2}\otimes a:a\in A\}\subset B(\ell_2\otimes  H)$
is a unital $w^*$-closed subalgebra of $B(\ell_2\otimes  H)$ 
with the $\A_1(1)$ property (see e.g., [Az, Section 2]).
Since $M_k(A^{(\infty)})=\bigl(M_k(A)\bigr)^{(\infty)}$, it follows that
$M_k\bigl(A^{(\infty)}\bigr)$ has the $\A_1(1)$ property for every $k\geq 1$. 

Applying Proposition 5.1 to $A^{(\infty)}$ and noticing that
$I_{\ell_2}\otimes A^{(\infty)}$ is canonically isomorphic to
$I_{\ell_2}\otimes A$, we obtain the following.

\proclaim{Corollary 5.3} Let $A\subset B(\H)$ be a unital, $w^*$-closed
subalgebra of $B(\H)$ and let $J\subset A$ be a $w^*$-closed $2$-sided
ideal.   Then there exists
a subspace $\E\subset \ell_2\otimes  \H$ such that the map
$\widehat{\Psi}:A/J\to B(\E)$ defined by $\widehat{\Psi}(a+J)
=P_\E (I_{\ell_2}\otimes a) {} {\vert_\E}$ is a 
completely isometric representation.
\endproclaim

\medbreak

An alternative proof for Corollary 2.10, i.e., 
  $\Phi: F^\infty/J\to P_{\N_J}F^\infty|_{\N_J}$ defined by $\Phi(f)= P_{\N_J}f|_{\N_J}$ is a completely isometric representation, can be obtained
  using Theorem 3.7 and Corollary 5.3 as follows. 

\demo{Alternative Proof of Corollary $2.10$}  From Corollary 5.3 (or from
Proposition 5.1 if we use that $F^\infty$ has property $\A_{\aleph_0}(1)$)
there exists a subspace $\E\subset \ell_2\otimes  \F^2$ such that the
map $\widehat{\Psi}:F^\infty/J\to B(\E)$, defined by
$\widehat{\Psi}(\eta)=P_\E\eta{}\vert_{\E}$, is a completely isometric
homomorphism.  Let $\varphi\in \E$ and notice that 
$\{I_{\ell_2}\otimes S_j : j\leq n\}$ 
satisfies $(2.1)$.  Then
$$
\sum_{|\alpha|=k}\|\widehat\Psi(S_\alpha+J)^*\varphi\|_2^2=
\sum_{|\alpha|=k}\|P_\E(I_{\ell_2}\otimes S_\alpha^*)\varphi\|_2^2\leq
\sum_{|\alpha|=k}\|(I_{\ell_2}\otimes S_\alpha^*)\varphi\|_2^2\to 0.$$
This shows that $[\widehat{\Psi}(S_1+J),\dots,
\widehat{\Psi}(S_n+J)]$ is $C_0$-contractive.

Notice that
for each $\varphi\in J$,
$\varphi \bigl(\widehat\Psi(S_1+J), \dots,\widehat\Psi(S_n+J) \bigr)=
\widehat\Psi(\varphi+J)=0$.  Then, from Theorem 3.7, there
exists a unital, completely contractive, $w^*$-continuous map 
$\Phi_K:B(\N_J)\to B(\E)$ satisfying
$\Phi_K(B_\alpha)= \widehat\Psi(S_\alpha+J)$ for every $\alpha\in\Fr_n^+$.  
Recall that $B_\alpha=\Phi(S_\alpha+J)$.  Hence, for each $\alpha\in\Fr_n^+$,
$\widehat\Psi(S_\alpha+J)=\Phi_K\circ\Phi(S_\alpha+J)$.  
Using the $w^*$-continuity of the three maps, we obtaing the
following commutative diagram 
$$\matrix  
F^\infty/J 
 & {\buildrel \hbox{$\widehat\Psi$}\over\longlongrightarrow}&B({\Cal E})\cr
\phantom{AAA}\Phi\searrow & & \nearrow \Phi_K\phantom{AAA}\cr
	 &B(\N_J)\cr
\endmatrix.$$
Since $\Phi_K$ and $\Phi$ are completely contractive,
and since $\widehat\Psi$ is completely isometric, we conclude that
$\Phi$ is completely isometric.
$\blacksquare$
\enddemo

Corollary 2.10 and the 
following simple lemma can be used to derive Theorem 2.4.  Thus,
we can prove this result without using the commutant lifting 
theorem of [Po3].  Notice that,
using a standard $w^*$-continuity argument, we can assume that the
$W_j$'s of Theorem 2.4 are $N\times N$ matrices.  We leave the
details to the reader.

\proclaim{Lemma 5.4} Let $\lambda_j=
(\lambda_{j1},\dots,\lambda_{jn})\in \Bbb B_n$,
$j=1,\dots,k$, be 
$k$ different points in $\Bbb B_n$. For each $i_0\in \{1,\dots,k\}$ 
there exists $\varphi_{i_0}\in F^\infty$ such that $\varphi_{i_0}(\lambda_{i_0})=1$ and
$\varphi_{i_0}(\lambda_j)=0$ whenever $i_0\not= j$.  Consequently, given
$W_1,\dots,W_k\in M_N$, there exists $\varphi\in M_N(F^\infty)$ such that
$\varphi(\lambda_j)=W_j$ for every $j=1,\dots,k$.
\endproclaim

\demo{Proof} Fix $i_0\in \{1,\dots,k\}$.
For each $j\not=i_0$ find $q\in \{1,\dots,n\}$ such that
$\lambda_{i_0q}\not=\lambda_{jq}$ and define $\theta_j=S_q-\lambda_{jq}I$.
Then $\theta_j(\lambda_{i_0})=\lambda_{i_0q}-\lambda_{jq}\not=0$ and
$\theta_j(\lambda_j)=0$.  
Let $\psi=\otimes_{j\not=i_0}\theta_j$.  Then
$\psi(\lambda_{i_0})\not=0$
and $\psi(\lambda_j)=0$ whenever $j\not=i_0$.  
Define $\varphi_{i_0}={1\over \psi(\lambda_{i_0})}\psi$.
If $W_1,\dots,W_k\in M_N$, then 
$\varphi=\sum_{i\leq k}W_i\otimes \varphi_i \in M_N(F^\infty)$
satisfies $\varphi(\lambda_i)=W_i$ for $i\leq k$.
$\blacksquare$
\enddemo

\enddocument